\documentclass[11pt]{article}
\usepackage{amssymb,amsmath,amsfonts,amscd,amsthm}

\usepackage{upgreek,color,url}
\usepackage{epsfig}
\usepackage{multirow}
\usepackage{fancybox}
\usepackage[textsize=tiny]{todonotes}
\usepackage{tikz}
\usetikzlibrary{backgrounds,matrix,fit}
\usetikzlibrary{calc}
\usetikzlibrary{positioning}
\usetikzlibrary{decorations.pathreplacing,angles,quotes}

\usepackage{blkarray}

\usepackage{algorithm}
\usepackage[noend]{algorithmic} 
\usepackage{authblk}
\usepackage[english]{babel} 
\usepackage{bm}
\usepackage{booktabs}    
\usepackage{hyperref,nameref,cleveref}
\usepackage{comment} 
\usepackage{empheq}
\usepackage[inline]{enumitem}
\usepackage{fancyhdr}
\usepackage[T1]{fontenc}    
\usepackage{fullpage}
\usepackage{gensymb}
\usepackage{graphicx}
\usepackage[utf8]{inputenc} 
\usepackage{latexsym}
\usepackage{lipsum}
\usepackage{listings} 
\usepackage{mathtools}
\usepackage{multicol,multirow}
\usepackage[square,sort,comma,numbers]{natbib}
\usepackage{nicefrac}  
\usepackage[parfill]{parskip}
\usepackage{setspace} 
\usepackage{subcaption}
\usepackage{times}
\usepackage{titlesec} 
\usepackage{thmtools}
\usepackage{upquote}
\usepackage{url}            
\usepackage{wrapfig}
\usepackage{xcolor}
\usepackage{xspace}
\usepackage{mystylefile}   
\usepackage{stackengine,trimclip}
\usepackage{colortbl}
\usepackage{pgfplots}
\pgfplotsset{compat=1.18}

\tikzset
{%
  pics/matrix/.style n args={6}{
    code={%
      \begin{scope}[y=-1cm,scale=0.5]
        \draw    (0,0) grid (#2,#1);
        \node at (0.5*#2,-0.5)   {#3};
        \node at (0.5*#2,#1+0.5) {#4};
        \node at (-0.5,0.5*#1)   {#5};
        \node at (#2+0.5,0.5*#1) {#6};
      \end{scope}     
    }},
}

    \theoremstyle{definition}

    \theoremstyle{definition}

    \theoremstyle{remark}
\newtheorem{remark}{Remark}
    \theoremstyle{remark}

    \theoremstyle{definition}

\newcommand{\E}{\mathbb{E}}
\newcommand{\F}{\mathbb{F}}

\newcommand{\Ab}{\mathbf{A}}

\newcommand{\Rb}{\mathbf{R}}

\newcommand{\Wb}{\mathbf{W}}

\newcommand{\Omegab}{\boldsymbol{\Omega}}

\newcommand*{\argmax}{\mathop{\mathrm{argmax}}}

\renewcommand{\diag}{\mathop{\mathrm{diag}}}

\newcommand{\nnz}{\mathop{\mathrm{nnz}}}

\newcommand{\nucnorm}[1]{{\left\vert\kern-0.25ex\left\vert\kern-0.25ex\left\vert #1 
    \right\vert\kern-0.25ex\right\vert\kern-0.25ex\right\vert}}

\newcommand{\wh}[1]{\widehat{#1}}
\newcommand{\wt}[1]{\widetilde{#1}}

\newcommand{\eg}{\emph{e.g.}\xspace}
\newcommand{\ie}{\emph{i.e.}\xspace}

\newcommand{\rbr}[1]{\left(#1\right)}
\newcommand{\sbr}[1]{\left[#1\right]}

\newcommand{\nbr}[1]{\left\|#1\right\|}
\newcommand{\abbr}[1]{\left|#1\right|}

\newcommand{\nnbr}[1]{{\left\vert\kern-0.25ex\left\vert\kern-0.25ex\left\vert #1 \right\vert\kern-0.25ex\right\vert\kern-0.25ex\right\vert}}

\newcommand{\bmat}[1]{\begin{bmatrix} #1 \end{bmatrix}}

\definecolor{commentcolor}{RGB}{110,154,155}   

\newcommand{\blue}[1]{\textcolor{blue}{#1}}



\newcommand{\algocom}[1]{\hfill \blue{$\triangleright$ #1}}

\title{Adaptive Parallelizable Algorithms for Interpolative Decompositions via Partially Pivoted LU}

\author{
    Katherine Pearce\thanks{
        Oden Institute, University of Texas at Austin (\texttt{katherine.pearce@austin.utexas.edu}).
    }
    \and
    Chao Chen\thanks{
        Department of Mathematics, North Carolina State University (\texttt{cchen49@ncsu.edu}).
    }
    \and
        Yijun Dong\thanks{
        Courant Institute, New York University (\texttt{yd1319@nyu.edu}).
    }
    \and
    Per-Gunnar Martinsson\thanks{
        Department of Mathematics \& Oden Institute, University of Texas at Austin (\texttt{pgm@oden.utexas.edu}).
    }
}
\date{}
\savestack\uppertriangle{%
\begin{tikzpicture}
  \node[outer sep=1in] (earmark) {};
  \filldraw[rounded corners, blue!40, opacity=0.20] (earmark.north west) -- (earmark.south east) -- 
    (earmark.north east)-- (earmark.north west);
  \end{tikzpicture}%
}

\begin{document}
\maketitle

\begin{center}
\begin{minipage}{135mm}
 \textbf{Abstract:} 
 Interpolative and CUR decompositions involve “natural bases” of row and column subsets, or skeletons, of a given matrix that approximately span its row and column spaces.  
 These low-rank decompositions preserve properties such as sparsity or non-negativity, 
 and are easily interpretable in the context of the original data. 
 For large-scale problems, randomized sketching to sample the row or column spaces with a random matrix can serve as an effective initial step in skeleton selection to reduce computational cost.
 A by now well established approach has been to extract a randomized sketch, followed by column-pivoted QR (CPQR) on the sketch matrix.
 This manuscript describes an alternative approach where CPQR is replaced by LU with partial pivoting (LUPP). 
 While LUPP by itself is not rank-revealing, it is demonstrated that when used in a randomized setting, LUPP not only reveals the numerical rank, but also allows the
 estimation of the residual error as the factorization is built. 
 The resulting algorithm is both adaptive and parallelizable, and attains much higher practical speed due to the lower communication requirements of LUPP over CPQR. 
 The method has been implemented for both CPUs and GPUs, and the resulting software has been made publicly available.
\end{minipage}
\end{center}


\section{Introduction}
\label{sec:intro}

Many problems in scientific computing involve matrices that are dense but ``data-sparse,'' or well-approximated by matrices of lower rank. 
In many applications, it is advantageous to work with a low-rank approximation that enables straightforward data interpretation in the same context as the input matrix. 
Interpolative and CUR decompositions \cite{mahoney2009, sorensen2014, woodruff2014optimalcur} are particularly well-suited for these applications because their approximate bases for the row and column spaces comprise actual rows and columns of the original matrix, often referred to as skeletons \cite{goreinov1995}. 
As such, these decompositions preserve not only individual matrix entries but also important properties such as sparsity or non-negativity. 

Applications in which the ID and CUR have been successfully utilized for data analysis include medical imaging and patient data analysis \cite{drineas2008relative}, information retrieval and recommendation systems \cite{drineas02}, and gene expression \cite{mahoney2009}. 
CUR and interpolative decompositions also arise more broadly in computational areas such as direct solvers for partial differential equations \cite{martinsson2019fast}.
We refer interested readers to \cite{mahoney2009} for an in-depth survey of ID and CUR decompositions and their applications.

There are many classical algorithms to select row or column skeletons greedily, often based on rank-revealing matrix decompositions like column-pivoted QR (CPQR) \cite[Section 5.4.1]{golub2013}, rank-revealing QR \cite{chan1987,ipsen94}, or strongly rank-revealing QR \cite[Algorithm~4]{gu1996}.
Despite their appealing empirical accuracy, greedy pivoting algorithms suffer from two known drawbacks. First, for those faster greedy pivoting methods like CPQR, there exist scarce adversarial inputs (\eg, Kahan-type matrices~\citep{kahan1966}) that significantly compromise the worst-case performance guarantee. Second, these greedy pivoting methods are inherently sequential as the skeleton selection in each step is based on the residual with respect to the previous skeleton selections.
For the first pitfall involving scarce adversarial inputs, replacing greedy pivoting with sampling~\citep{frieze2004fast,belabbas2009spectral,cohen2015uniform,derezinski2021determinantal,mahoney2009cur} or random pivoting~\citep{deshpande2006matrix,deshpande2006adaptive,chen2022randomly,dong2023robust} serves as a simple but effective remedy.
Meanwhile, the latter setback has been mitigated by the development of randomized sketching for dimension reduction, in which random matrices are applied to the input matrix to sample its row or column space \cite{halko2011, liberty20167, martinsson2020, chen2020, dong2022}. 
Randomized sketching also lends itself naturally to parallelization for high-performance computing \cite{dongarra20, martquin19}.

The challenges of developing a high-performance parallelized version of QR, however, are inherent in its algorithmic implementations \cite{Henry1996ParallelizingTQ, demgrig12, martquin19}. 
The performance of column-pivoted Householder QR, for instance, is capped by the number of operations (approximately half; see \cite{quintanaorti98}) that can be cast in terms of Level-3 Basic Linear Algebra Subprograms
(BLAS) \cite{dongarra90, dongarra91}, high-performance matrix-matrix multiplications.
Since the GPU efficiently handles compute-intensive tasks, the memory-bound Level-2 and Level-1 BLAS operations prevalent in classical CPQR algorithms impede maximally-leveraged GPU computations. 
Issues that also arise in parallelizing QR algorithms include inefficient storage with unevenly distributed workloads among the processors as well as communication bottlenecks \cite{Henry1996ParallelizingTQ}. 
The difficulties in designing parallel pivoting algorithms that maximize GPU computing resources can be circumvented if CPQR is replaced by partially pivoted LU (LUPP). However, LUPP does not reveal the numerical rank of a matrix (CPQR does that with strong theoretical guarantees \cite{chan1987, ipsen94, gu1996}) because it can fail to pivot so that smallest-magnitude entries are confined to the bottom right of the upper triangular matrix $\mU$ (see (\ref{eq:lupp}) of Section~\ref{sec:lupp}) \cite{peters1975, pan2000rrlu}. 

\begin{figure}
\centering
%
\centering
\includegraphics[width=0.35\textwidth]{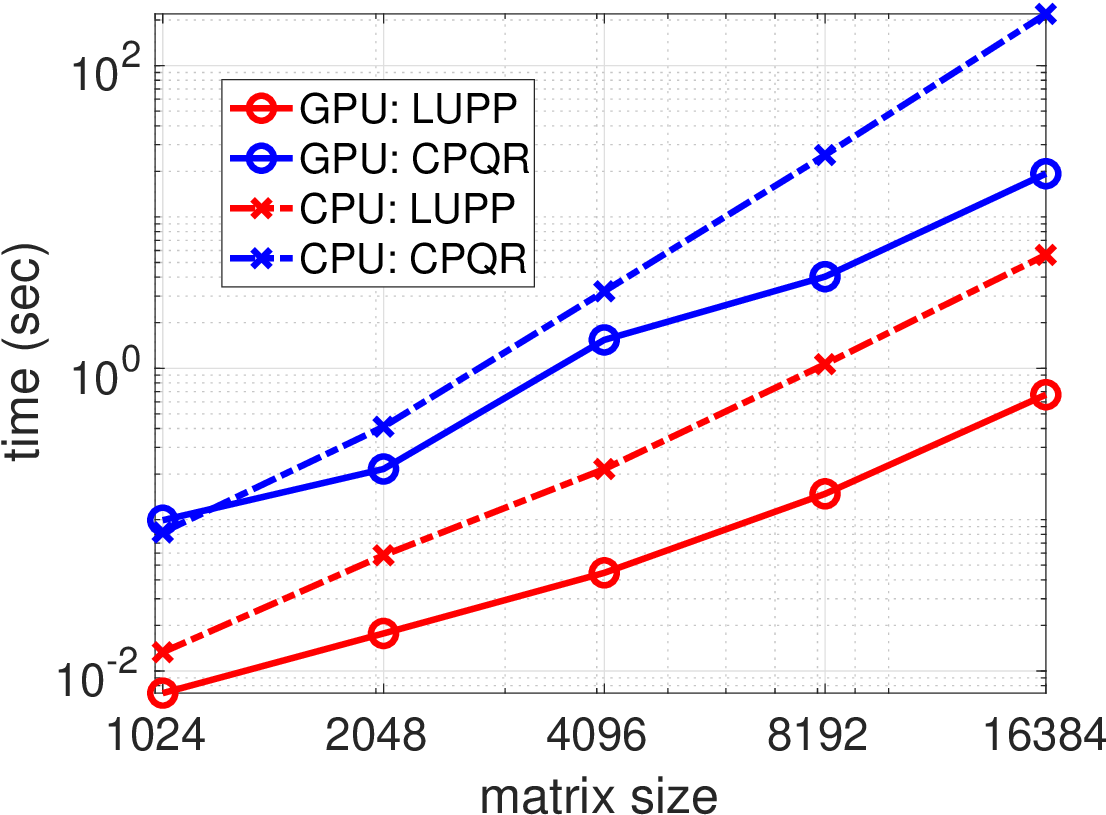}
\caption{Running time of parallel LUPP and parallel CPQR on an NVIDIA V100 GPU and on two Intel Xeon Gold 6254 CPUs.}
\label{f:luqr}
\end{figure}

It has recently been shown that LUPP in conjunction with randomized sketching, e.g. \cite[Algorithm 2]{dong2022}, exhibits rank-revealing properties comparable to CPQR \cite[Section 5]{dong2022}. 
Additionally, LUPP is more computationally efficient than CPQR and reaps more benefits from parallelization and GPU-acceleration \cite{grigdem08, demgrig12, FAVERGE2015}; see also
Figure~\ref{f:luqr}.

However, if the target rank, or embedding dimension for the randomized sketch, is unknown a priori, the lack of rank-revealing guarantees for LUPP have obfuscated an incremental blocked version based on a certificate of accuracy like \cite[Algorithm 13]{martinsson2020}.
 Computing orthogonal projections to estimate the low-rank approximation error in an incremental blocked LUPP algorithm is tantamount to performing both LU and QR, undermining the performance improvements from using LUPP.

With randomized sketching as an initial step, we show that LUPP can be modified into an adaptive, parallelized skeleton selection method with cost-efficient estimates for the certificate of accuracy.
In Section~\ref{sec:bg}, we provide an overview of relevant numerical linear algebra material, including interpolative and CUR decompositions and common algorithms for skeleton selection.
Our contributions include a technique to estimate the low-rank approximation error associated with these decompositions in Section~\ref{sec:errest}, which we use as an error estimate in our adaptive randomized algorithm presented in Section~\ref{sec:blkLUPP}. 
Accuracy comparisons between skeletons selected by our algorithm and skeletons selected by randomized CPQR are given in Section~\ref{sec:numexp}, as well as numerical evidence which supports that the asymptotic bound on the growth factor of Haar-distributed random matrices found in \cite{Higham21} holds for a larger class of random matrices.
We also provide comparisons of our algorithm's performance to state-of-the-art algorithms with both CPU and GPU implementations. 
We conclude with a summary of our investigation and potential avenues for future work in Section~\ref{sec:concl}.

\paragraph{Summary of Contributions}
\begin{enumerate}
    \item We introduce an error estimation technique that enables an adaptive, randomized approach to determine skeletons for CUR and interpolative decompositions.
    \item We present an adaptive randomized algorithm based on our error estimate with LUPP to find CUR and IDs that can be parallelized to leverage GPU computations.
    \item We demonstrate that our algorithm yields an ID of comparable accuracy to randomized algorithms based on CPQR.
    \item We illustrate the improved computational efficiency of our algorithm over existing algorithms based on CPQR, with performance comparisons for both CPU and GPU implementations.
    \item We provide software for CPU and GPU implementations.
    \item We propose a blockwise adaptive sketching scheme for randomized LUPP with a set of efficient posterior error estimates. This enables automatic target rank detection for any given error tolerance while leveraging the desired parallelizability of both sketching and LUPP.
\end{enumerate}

\section{Preliminaries}
\label{sec:bg}

We briefly review relevant concepts from numerical linear algebra used throughout the text and summarize  notation in Section~\ref{sec:not}. We then review ID and CUR decompositions in Section~\ref{sec:ID}.
In Section~\ref{sec:piv}, we demonstrate how such decompositions may be obtained from row or column pivoting algorithms.
Lastly, in Section~\ref{sec:emb}, we give a brief summary of randomized sketching, the step that renders a ``rank-revealing'' LUPP. 

\subsection{Notation}
\label{sec:not}

We use the Euclidean norm to measure vectors $\xb \in \F^{n}$ over a field $\F$.
Matrices $\ma \in \F^{m \times n}$ are measured in the Frobenius norm $\| \ma\|_F = (\sum_{i,j} |\ma(i,j)|^2)^{1/2}$.
We adopt the notation of Golub and Van Loan~\cite{golub2013} to refer to submatrices: 
for subsets of row and column index vectors $I = \bmat{i_1 & i_2 & \ldots, i_k}$ and $J = \bmat{j_1 & j_2 & \ldots, j_{\ell}}$, let $\ma(I,J)$ denote the $k \times \ell$ matrix
\begin{align*}
    \ma(I,J) = \bmat{\ma(i_1,j_1) & \ma(i_1,j_2) & \cdots & \ma(i_1,j_{\ell} \\ \ma(i_2,j_1) & \ma(i_2,j_2) & \cdots & \ma(i_2,j_{\ell} \\ \vdots & \vdots & & \vdots \\ \ma(i_k,j_1) & \ma(i_k,j_2) & \cdots & \ma(i_k,j_{\ell}}.
\end{align*}
The matrix $\ma(I,:)$ is shorthand for the matrix $\ma(I,\bmat{1 & \ldots & n})$, and we define $\ma(:,J)$ analogously.
We let $\ma^*$ denote the transpose (or Hermitian transpose) of $\ma$, and we say a matrix $\mU$ is \textit{orthonormal} if its columns are orthonormal, so that $\mU^*\mU = \mI$.
The Moore-Penrose pseudoinverse of $\ma$ is denoted by $\ma^\dag$.

\subsubsection{The QR factorization} Every matrix $\ma \in \F^{m \times n}$ admits a \textit{QR factorization} given by
\begin{align}
    \label{eq:QR}
    \begin{array}{ccccc}
       \ma  & \mP & = & \mQ & \mR ,\\
       m \times n  & n \times n &  & m \times r &  r \times n
    \end{array}
\end{align}
where $r = \min(m,n)$, $\mQ$ is orthonormal, $\mR$ is upper triangular, and $\mP$ is a permutation matrix.
In practice, we represent $\mP$ by a permutation vector $J$ of column indices, so that $\mP = \mI(:,J)$, where $\mI$ is the $n \times n$ identity matrix.
Then the factorization (\ref{eq:QR}) can be expressed more succinctly as
\begin{align*}
    \begin{array}{cccc}
       \ma(:,J)  & = & \mQ & \mR .\\
       m \times n  &  & m \times r &  r \times n
    \end{array}
\end{align*}
Using a greedy algorithm that builds the QR factorization incrementally
(e.g. Gram-Schmidt), we can halt the factorization after computation of the first $k$ terms to obtain a ``partial'' QR factorization:
\begin{align}
\label{eq:qrtrunc}
    \begin{array}{cccc}
       \ma(:,J)  & = & \mQ_k & \mR_k .\\
       m \times n  &  & m \times k &  k \times n
    \end{array}
\end{align}
We define the functions 
\begin{align}
\label{eq:qrfun}
    \bmat{\mQ, \mR, J} = \texttt{qr}(\ma) \ \ \textup{and} \ \ \bmat{\mQ, \mR, J} = \texttt{qr}(\ma, k),
\end{align}
that produce either the full or partial QR factorizations of $\ma$, respectively.

\subsubsection{The singular value decomposition (SVD)} Every matrix $\ma \in \F^{m \times n}$ also admits a \textit{singular value decomposition}, given by
\begin{align}
    \label{eq:svd}
    \begin{array}{ccccc}
       \ma  & = & \mU & \mathbf{\Sigma} & \mV^*,\\
       m \times n  &  & m \times r &  r \times r & r \times n
    \end{array}
\end{align}
where $r = \min(m,n)$, matrices $\mU$ and $\mV$ are orthonormal, and $\mathbf{\Sigma}$ is diagonal.
The columns $\{\uu_{i}\}_{i=1}^{r}$ and $\{\vv_{i}\}_{i=1}^{r}$ of $\mU$ and $\mV$, respectively, are called the left and right singular vectors of $\ma$.
The diagonal elements $\{ \sigma_i\}_{i=1}^{r}$ of $\mathbf{\Sigma}$ are the singular values of $\ma$, ordered so that $\sigma_1 \geq \sigma_2 \geq \cdots \geq \sigma_r \geq 0$.
We can again truncate the factorization after the first $k$ terms, letting $\ma_k$ denote $\sum_{i=1}^k \sigma_i \uu_i \vv_j^*$.
By the Eckart-Young theorem \cite{eckart1936}, the singular values yield the smallest errors incurred in a rank-$k$ approximation of $\ma$:
\begin{align}
    \|\ma - \ma_k \|_2 = \sigma_{k+1} \ \ \textup{and} \  \ \|\ma - \ma_k\|_F = \left ( \sum_{j = k+1}^{r} \sigma_j^2 \right )^{1/2}.
\end{align}

\subsection{Interpolative and CUR Decompositions}
\label{sec:ID}

Interpolative decompositions (ID) and the CUR decomposition are formed using actual rows and columns of a matrix $\ma$, so-called ``natural bases'' approximately spanning its row and column spaces. 
Remaining rows or columns of $\ma$ are expressed in terms of the selected ones.
Low-rank matrix approximations that utilize these decompositions preserve properties such as sparsity or non-negativity.
They are also memory-efficient in the sense that only the selected row or column indices, rather than matrix entries, need to be stored. 
We briefly describe each decomposition before discussing approximation error bounds.

\subsubsection{Definitions}
\label{sec:IDdef}

We assume throughout this section $\ma \in \F^{m \times n}$ has rank at least $k$.

\paragraph{Row ID} Given $k$ linearly independent rows of $\ma$, indexed by $I$, the row ID of $\ma$ has the form
\begin{align}
    \label{eq:rowID}
    \begin{array}{cccc}
       \ma  & = & \mW & \ma(I,:) ,\\
       m \times n  &  & m \times k &  k \times n
    \end{array}
\end{align}
where the \textit{row interpolation matrix} $\mW$ contains the $k \times k$ identity matrix as a submatrix, corresponding to the rows of $\ma$ indexed by $I$, also known as skeletons.
Moreover, every entry of $\mW$ is bounded in magnitude by 1, though in practice this bound may be slightly exceeded. 
For the corresponding error bounds, if $\mR = \ma(I,:)$, it is useful to express (\ref{eq:rowID}) in terms of the orthogonal projector $\mR^{\dag}\mR$ onto the span of the selected rows:
\begin{align}
    \label{eq:rowID_orth}
    \begin{array}{cccc}
       \ma  & = & \ma \mR^\dag & \mR.\\
       m \times n  &  & m \times k &  k \times n
    \end{array}
\end{align}

Thus, we are interested in the row ID error 
\begin{align}
\label{eq:IDerr_row}
    \| \ma - \ma \mR^\dag \mR \|_F.
\end{align}
Because $\mR$ may have a large condition number, the row ID error can be computed in a numerically stable manner with orthonormalization:
\begin{align}
    \label{eq:sID_row}
    \| \ma - \ma \mQ \mQ^* \|_F,
\end{align}
where $\mQ$ is the output of \texttt{qr}
applied to the rows of $\mR = \ma(I,:)$.
Without loss of generality, we will compute row IDs to illustrate our adaptive algorithm, but we also include short summaries here of the analogously computed column ID, two-sided ID, and CUR decomposition for completeness. 

\paragraph{Column ID} Given $k$ linearly independent columns of $\ma$ indexed by $J$, the column ID of $\ma$ is
\begin{align*}
    \begin{array}{cccc}
       \ma  & = & \ma(:,J) & \mX ,\\
       m \times n  &  & m \times k &  k \times n
    \end{array}
\end{align*}
where the \textit{column interpolation matrix} $\mX$ contains the $k \times k$ identity matrix as a submatrix, corresponding to the column skeletons  indexed by $J$; its entries are also bounded in magnitude by 1. 
Letting $\mC = \ma(:,J)$, we can re-write the column ID in terms of the orthogonal projector $\mC \mC^\dag$ onto the span of the selected columns:
\begin{align*}
    \begin{array}{cccc}
       \ma  & = & \mC & \mC^\dag \ma.\\
       m \times n  &  & m \times k &  k \times n
    \end{array}
\end{align*}
The ID approximation error $\|\ma - \mC \mC^\dag \ma \|_F$ can be computed in a numerically stable method as above by orthonormalizing $\mC$. 

\paragraph{Two-sided ID} Given $k$ linearly independent rows of $\ma$ indexed by $I$, and $k$ linearly independent columns of $\ma$ indexed by $J$, define $\ms = \ma(I,J)$ as the submatrix of skeletons.
The two-sided ID is
\begin{align*}
    \begin{array}{ccccc}
       \ma  & = & \mC \ms^{-1} & \ms & \mC^{\dag} \ma, \\
       m \times n  &  & m \times k & k \times k &  k \times n
    \end{array}
\end{align*}
for $\mC = \ma(:,J)$.
In exact arithmetic, $\ms^{-1} \ms = \mI$, so the two-sided ID is equivalent to the column ID. 

\paragraph{CUR} Given  $k$ row skeletons indexed by $I$ and $k$ column skeletons indexed by $J$, the CUR decomposition of $\ma$ is given by
\begin{align*}
    \begin{array}{ccccc}
       \ma  & = & \mC  & \mC^\dag \ma \mR^\dag & \mR, \\
       m \times n  &  & m \times k & k \times k &  k \times n
    \end{array}
\end{align*}
While CUR resembles the two-sided ID, their approximation errors and conditioning differ.
In general, construction of the CUR decomposition is more ill-conditioned than the two-sided ID, due to the typically large condition numbers of $\mC$, $\mR$, and $\mC^\dag \ma \mR^\dag$.
In the two-sided ID, only the computation of $\ms$ may potentially be ill-conditioned;
evaluations of $\mC \ms^{-1}$ and $\mC^\dag \ma$ are well-conditioned, and explicit computation of direct inverses can be avoided in practice~\citep{martinsson2019randomized}.
We refer interested readers to \cite{sorensen2014,voronin2017,anderson2015spectral,martinsson2019randomized} for more details on the conditioning of these decompositions.

\subsection{Pivoting Algorithms for Matrix Decompositions}
\label{sec:piv}

In this section, we review two algorithms utilized to determine row and column skeletons. 
We assume throughout that $\ma \in \F^{m \times n}$ has rank $k$, and $r = \min(m,n) = n$ for ease of explanation.

\subsubsection{QR with Column Pivoting (CPQR)}
\label{sec:cpqr}

The CPQR algorithm underlies the \texttt{qr} function from Section~\ref{sec:bg}.
We recall \texttt{qr} yields the factorization $\ma \mP = \mQ \mR$, which is
achieved through a sequence of permutations and Householder reflections that rank-1 update the ``active submatrix'' of $\ma$ \cite{householder1958hhqr}.
Formally, let $\ma^{(t)}$ represent the resulting matrix after $t$ steps of pivoting and updating for $t = 0,\ldots,n-2$ (with $\ma^{(0)} = \ma$).
At the $(t+1)$-st step, CPQR identifies the column pivot $j_{t+1}$ having maximal Euclidean norm in the active submatrix,
\begin{align*}
    j_{t+1} = \argmax_{t+1 \leq j \leq n} \left \| \ma^{(t)}(t+1:m, j) \right \|_2.
\end{align*}

Partitioning so that 
\begin{align}
    \label{eq:cpqr}
    \begin{array}{ccccccc}
            &     &  &  k \hspace{3mm} n-k  &  k \hspace{3mm} n-k &  \\
       \ma  & \mP & =  & \underbrace{\bmat{ \mQ_1  & \mQ_2 }}_{\mQ} & \underbrace{\bmat{\mR_{11} & \mR_{12} \\ \mzero & \mR_{22}}}_{\mR} &  \begin{array}{c}  k \\ n - k \end{array} ,
    \end{array}
\end{align}
and using Theorem 7.2 of \cite{gu1996}, the entries $ \left | \left (\mR_{11}^{-1} \mR_{12} \right )_{ij} \right| \leq 2^{k-i}$ for $1 \leq i \leq k$ and $1 \leq j \leq n-k$; additionally, the relevant singular value gaps satisfy
\begin{align*}
    \sigma_k(\mR_{11}) \geq \frac{\sigma_k(\ma)}{\sqrt{n-k} \, 2^k} \ \ \ \textup{and} \ \ \ \sigma_1(\mR_{22}) \leq \sqrt{n-k} \, 2^k \sigma_{k+1}(\ma). 
\end{align*}
While there are adversarial cases that achieve these exponential bounds like Kahan-type matrices~\citep{kahan1966, peters1975}, they are rarely encountered in practice~\cite{trefethen1990}. 




\subsubsection{LU with Partial Pivoting (LUPP)}
\label{sec:lupp}

Following the notation of \cite{golub2013}, the LU algorithm with row-wise partial pivoting (LUPP) yields
\begin{align}
    \label{eq:lupp}
    \begin{array}{cccccc}
              &   &     &   k \hspace{3mm} n-k  &  k \hspace{3mm} n-k   \\
       \mP  \ma & = & \begin{array}{c} k \\ m-k \end{array} & \underbrace{\bmat{ \mL_{11}  & \mzero \\ \mL_{21} & \mL_{22} }}_{\mL} & \underbrace{\bmat{\mU_{11} & \mU_{12} \\ \mzero & \mU_{22}}}_{\mU},
    \end{array}
\end{align}
where $\mL$ is lower triangular with main diagonal entries equal to 1 and off-diagonal entries bounded in magnitude by 1, and $\mU$ is upper triangular.

The pivoting strategy to construct this factorization updates the active submatrix via Schur complements (e.g. \cite[Algorithm 3.2.1]{golub2013} \cite[Algorithm 21.1]{trefethen1997}).
Let $\ma^{(t)}$ be the resulting matrix after $t$ steps of pivoting and updating for $t = 0, 1, \ldots, n-2$. 
At the ($t+1$)-st step of LUPP, the largest-magnitude column entry of the active submatrix is selected as the next pivot entry:
\begin{align*}
    i_{t+1} = \argmax_{t+1 \leq i \leq m} |\ma^{(t)}(i,t+1)|.
\end{align*}
The active submatrix is then updated: defining vectors $\xr = t+2:m$ and $\xc = t+2:n$,
\begin{align}
    \label{eq:lupps}
    \ma^{(t+1)}(\xr, \xc) = \ma^{(t)}(\xr, \xc) - \frac{\ma^{(t)}(\xr,t+1) \ma^{(t)}(t+1, \xc)}{\ma^{(t)}(t+1,t+1)},
\end{align}
where $\ma^{(t)}(t+1,t+1)$ has been updated with $\ma^{(t)}(i_{t+1},t+1)$ following row permutation.

As in CPQR, the LUPP algorithm leads to an exponential upper bound on the entries of $\mU$, usually expressed in terms of the growth factor $\rho = \max_{i,j} |\mU(i,j)|/ \max_{i,j}|\ma(i,j)| \leq 2^{n-1}$. 
This upper bound is tight for certain matrices (e.g. Kahan-type matrices~\citep{kahan1966, peters1975, sorensen2014}), but again these types of adversarial inputs are rare in practice~\citep{trefethen1990}. 
While LUPP can fail for rank-deficient matrices and is notably not rank-revealing, \cite{trefethen1990} hypothesizes that the average-case growth factor grows sublinearly with respect to the size of the input matrix. 
Despite potential drawbacks, LUPP exhibits superior computational efficiency and parallelizability over CPQR. 
Additionally, randomization techniques like sketching stabilize LUPP in practice~\cite{trefethen1990,pan2015rand, pan2017num, dong2022}. 


\subsection{Matrix Sketching}
\label{sec:emb}

Given a matrix $\ma \in \F^{m \times n}$ of rank $k \ll \min(m,n)$, and a randomized linear embedding $\mathbf{\Omega} \in \F^{n \times \ell}$ with $k < \ell \ll \min(m,n)$, we define a \textit{randomized column sketch} of $\ma$ as the matrix $\mY = \ma \mathbf{\Omega} \in \F^{m \times \ell}$. 
Similarly, given $\mathbf{\Gamma} \in \F^{m \times \ell}$, we define a \textit{randomized row sketch} of $\ma$ to be $\mX = \ma^* \mathbf{\Gamma} \in \F^{n \times \ell}$. 
The process of forming $\ma \mathbf{\Omega}$ or $\ma^* \mathbf{\Gamma}$ is called sketching \cite{alonmat99, martinsson2020}, a well-established technique for dimension reduction \cite{jlind84, dasgupta2003jlproof, halko2011, martinsson2019randomized, demmel2019improved} used in a variety of applications \cite{linrab95, indyk1998, alongib02, rokhlin2008, levitt2022randomized}. 
In particular, it has facilitated computations where classical techniques in numerical linear algebra have proven too expensive or intractable for modern computing architectures \cite{Drineas06a, Drineas06b, Drineas06c, GhaWood16a, GhaPhi16b}. 

\section{Error Estimation with Randomized Sketches}
\label{sec:errest}

In this section, we demonstrate how randomized sketching enables efficient and reliable computation of ID and CUR factorizations and provides an estimate of their low-rank approximation errors. 

\subsection{Randomized Sketching for ID and CUR Factorizations}
\label{sec:rs_ID}

Rather than applying column or row skeleton selection algorithms directly to the input matrix $\ma \in \F^{m \times n}$, we first sketch with random matrices $\mathbf{\Omega} \in \F^{n \times \ell} $ or $\mathbf{\Gamma} \in \F^{m \times \ell}$ as in Section~\ref{sec:emb} to reduce the problem size.
Once we have the random sample $\mY = \ma \mathbf{\Omega}$ or $\mX = \ma^* \mathbf{\Gamma}$, we apply skeleton selection algorithms to the sample matrix.

To illustrate, suppose we have the sample matrix $\mY = \ma \mathbf{\Omega} \in \F^{m \times \ell}$.
To find a CUR factorization of $\ma$, we apply LUPP to $\mY$ to determine the desired row indices $I = \bmat{i_1, i_2, \ldots, i_{\ell}}$.
These indices directly correspond to row indices of the original matrix $\ma$, and we form the submatrix $\ma(I,:) \in \F^{\ell \times n}$.
Since $\ell \ll \min(m,n)$ by assumption, it is relatively inexpensive to determine column skeletons $J = \bmat{j_1, j_2, \ldots, j_{\ell}}$ by applying LUPP to  $\ma(I,:)^*$. 
The ID errors in Section~\ref{sec:ID} may then be computed with $I$ and $J$ since the formulas are independent of the method of skeleton selection.
However, we aim to use a random sample of the ID error in a more efficient method of error estimation, which can also be used in an adaptive algorithm.

To maximize computational efficiency, we opt for LUPP in the procedure illustrated above, referred to as randomized LUPP or \texttt{randLUPP} \cite[Algorithm~2]{dong2022}.
We refer to a randomized skeleton selection algorithm based on CPQR as \texttt{randCPQR} (see \cite[Algorithm~1]{dong2022}, \cite{voronin2017, sorensen2014}).
With the same goal of empirical efficiency, we describe a method of approximation error estimation employed in our adaptive version that leverages the LUPP factorizations of the sample matrix $\mY$ from the previous iterations. 

\subsection{Error Estimation with the Schur Complement}
\label{sec:errsch}

We revisit the factorization in (\ref{eq:lupp}), now with the sample matrix $\mY = \ma \mathbf{\Omega} \in \F^{m \times k}$ for some target rank $k$.
We note that we use $k$ in this section, rather than $\ell$ as we did in introducing randomized sketching, for consistency with our notation in the adaptive version of the algorithm in later sections (see also Remark~\ref{sec:remark2}).

For a row ID computation, we consider the LUPP factorization of $\mY$  given by
\begin{align}
    \label{eq:luppy}
    \begin{array}{ccccccccc}
    \mP & \mY & = &  \mL & \mU & = & \begin{array}{c} k \\ m - k \end{array} & \bmat{\mL_1 \\ \mL_2} & \mU \\
    m \times m & m \times k &  & m \times k & k \times k & &  & 
    \end{array}
\end{align}

The permutation matrix $\mP$ identifies  $k$ row pivots of $\mY$, denoted by the row index vector $I = \bmat{i_1, i_2, \ldots, i_{k}} \subset \bmat{1,\ldots,m}$.
Let $\mW$ be the row interpolation matrix containing the $k \times k$ identity matrix as a submatrix, corresponding to the row indices $I$, so that
\begin{align}
    \label{eq:inty}
    \mP \mW = \bmat{\mI_{k \times k} \\ \mL_{2} \mL_{1}^{-1}}.
\end{align}
Define $\mR = \ma(I,:) = (\mP \ma)(1:k,:)$, and consider the product 
\begin{align*}
    \mP \mW \mR = \bmat{(\mP \mY)(1:k,:) \\ (\mL_{2} \mL_{1}^{-1})((\mP \mY)(1:k,:))}
\end{align*}

To estimate the approximation error associated with the row ID (\ref{eq:rowID}), we use the ``certificate of accuracy'' idea introduced in 
\cite[Sec.~12]{martinsson2020}, which is based on a randomized technique for estimating matrix norms. The quantity that we seek to bound is 
$$
    \|\ma - \ma \mR^\dag \mR \|_F = \|\ma - \mW \mR\|_F \\
    = \|\mP \ma - \mP \mW \mR \|_F .
$$
We now draw an independent random sample  $\wt{\mathbf{\Omega}} \in \F^{n \times b}$ for some $b < k$ from the Gaussian distribution with standard deviation $1/b$. Then for any matrix $\mtx{X}$
of size $m\times n$, the expectation of $\|\mtx{X}\wt{\mathbf{\Omega}}\|_{\rm F}^{2}$
is $\|\mtx{X}\|_{\rm F}^{2}$, with a modest variance. In consequence, we can estimate the norm of $\mP \ma - \mP \mW \mR$ by forming the independent sample $\mZ = \ma \wt{\mathbf{\Omega}} \in \F^{m \times b}$, since 
\begin{multline*}
    \|\ma - \ma \mR^\dag \mR \|_F 
    = \|\mP \ma - \mP \mW \mR \|_F 
    \approx \|\mP \ma \wt{\mathbf{\Omega}} - \mP \mW \mR \wt{\mathbf{\Omega}} \|_F \\
    = \|(\mP \mZ)(k+1:m, :) -  (\mL_{2} \mL_{1}^{-1})((\mP \mZ)(1:k,:)) \|_{\rm F}
\end{multline*}
The Schur complement used for the error estimate is then 
\begin{align}
\begin{split}
\label{eq:errSch}
   E_{Schur} &= \|(\mP \mZ)(k+1:m,:) - \mL_{2} \mL_{1}^{-1}((\mP \mZ)(1:k,:)) \|_F.
\end{split}
\end{align}

\remark{ \label{sec:remark1} We have observed in practice that a ``residual'' upper triangular matrix $\mU_{r} \in \real^{p \times p}$ for some small oversampling parameter $p$ (e.g.~$5 \leq p \leq 10$), given by
\begin{align*}
    \mP [ \underbrace{\mY}_{m \times k} \ \underbrace{\mY_r}_{m \times p}] &= \mL \mU \\
    &= \bmat{\mL_{1} & \\ \mL_2 & \mL_r} \bmat{\mU_1 & \mU_2 \\  & \mU_r}, 
\end{align*}
can also serve as a reliable proxy for the error estimate, where $\mY_r = \ma \mOmega_r \in \real^{m \times p}$ is an independent random sample set aside for this estimation, and $\mOmega_r$ is Gaussian with standard deviation $1$.
Then as above,
\begin{align*}
    \|\ma - \ma \mR^\dag \mR \|_F 
    = \|\mP \ma - \mP \mW \mR \|_F 
    \approx \frac{1}{\sqrt{p}} \|\mP \ma \mOmega_r - \mP \mW \mR \mOmega_r \|_F \\
    = \frac{1}{\sqrt{p}} \|(\mP \mY_r)(k+1:m, :) -  (\mL_{2} \mL_{1}^{-1})((\mP \mY_r)(1:k,:)) \|_{\rm F} \\
    =  \frac{1}{\sqrt{p}} \| (\mL_2 \mU_2 + \mL_r \mU_r) - (\mL_2 \mL_1^{-1})(\mL_1 \mU_2)  \|_F 
    = \frac{1}{\sqrt{p}} \| \mL_r \mU_r \|_F \\
    \leq \frac{1}{\sqrt{p}} \  \| \mL_r \|_2 \| \mU_r \|_F \leq \frac{1}{\sqrt{p}} \sqrt{p(m-k)} \ \|\mU_r \|_F \\
    = \sqrt{m-k} \ \| \mU_r \|_F.
\end{align*}
However, this upper bound is very pessimistic and not informative on its own. 
Curiously, we have found that it can be tightened to give a reliable error estimate
\begin{align}
    \label{eq:normUest}
    \| \ma - \ma \mR^\dag \mR  \|_F \lesssim \frac{4 \log k}{k} \sqrt{m-k} \ \|\mU_r\|_F,
\end{align}
using the reciprocal of the asymptotic upper bound $\frac{k}{4 \log k}$ on the LU growth factor of Haar-distributed random matrices introduced in \cite{Higham21}.
More curiously, this error estimate holds for other random matrix distributions including sparse sign matrices \cite{meng2013, nelson2013, woodruff2015} and subsampled randomized trigonometric transforms \cite{halko2011, rokhlin2008, tropp2011} (with a factor of $\sqrt{p}$).
It is worth noting that the bound can be further improved using the largest-magnitude element of $\mU_r$,
\begin{align}
    \label{eq:maxUest}
    \| \ma - \ma_k \| \lesssim \frac{4 \log k}{k} \sqrt{m-k} \ \max_{1 \leq i,j \leq p} |\mU_r(i,j)|,
\end{align}
which is more robust to oversampling than (\ref{eq:normUest}), but may provide an underestimate of the error; 
note $\max_{i,j}|\mU_r(i,j)|$ is a lower bound on $\max_{i,j}|\mY_r^{(\ell)}(i,j)|$, where $\mY_r^{(\ell)}$ is the active submatrix in the $\ell$-th step for $\ell = 1,\ldots,p-1$ in Gaussian elimination with $\mY_r$.
This issue can be circumvented by computing an a posteriori certificate of accuracy as in \cite{martinsson2020} using the selected skeleton indices upon termination of the algorithm.
Given the numerical evidence in Section~\ref{sec:acccomp} and Appendix~\ref{sec:appendix} for a variety of input matrices and random matrix distributions, we hypothesize that the upper bound on the growth factor found in \cite{Higham21} holds more generally than for the Haar distribution, also supporting the recent conjecture in \cite{dong2023robust} and the analysis of sub-Gaussian distributions in \cite{saibaba2023randomized}.
}


\section{Adaptive Randomized LUPP for ID/CUR Factorizations}
\label{sec:blkLUPP}

In the previous sections, we assumed \textit{a priori} knowledge of a randomized embedding dimension $k$. 
Here we present an adaptive version of \texttt{randLUPP} that relies on the error estimate (\ref{eq:errSch}) that can also determine an appropriate value of $k$ for a rank-$k$ ID or CUR factorization.

More concretely, in the adaptive version of randomized LUPP, \texttt{randLUPPadap} (Algorithm~\ref{alg:IRMS}), we seek a cost-effective dimension $k_t$ after $t$ iterations that yields an approximate basis of skeletons for the row or column space of $\ma \in \F^{m \times n}$.
For the row ID, for instance, upon termination of the algorithm after $t$ iterations, we have an integer $k_t$, a subset $I^{(t)}_s$ of $k_t$ row skeletons, and an interpolation matrix $\mW^{(t)}$ satisfying
\begin{align*}
    \|\ma - \mW^{(t)}\mR^{(t)} \|_F \lesssim \tau,
\end{align*}
for a prescribed error tolerance $\tau$, where $\mR^{(t)} = \ma(I_s^{(t)},:)$. 
The procedure to determine column skeletons  is analogous, using $\ma^* \mathbf{\Gamma}^*$. 

Recall that \texttt{randLUPP} determines a row ID by forming a sample matrix $\mY = \ma \mathbf{\Omega} \in \F^{m \times k}$ and applying LU with partial pivoting to $\mY$ to find row skeletons $I_s = \bmat{i_1,\ldots, i_{k}}$.
In \texttt{randLUPPadap}, we form random samples $\mY^{(t)} = \ma \mathbf{\Omega}^{(t)} \in \F^{m \times b}$ for $t = 0,1,2,\ldots$ with independent random matrices $\mathbf{\Omega}^{(t)} \in \F^{n \times b}$ using a prescribed block size $b$.
We compute a partially pivoted LU factorization of each random sample $\mY^{(t)}$ to arrive at a partially pivoted LU factorization of the full sample matrix $\mY = \bmat{\mY^{(0)}, \mY^{(1)},\ldots,\mY^{(t)}} \in \F^{m \times (t+1)b}$. 
The algorithm terminates when the factorization of the full sample matrix in that iteration yields a row skeleton vector $I^{(t)}_s$ for which $E^{(t)}_{Schur}$ in (\ref{eq:errSch}) satisfies the given error tolerance. 

We describe the algorithm's first iteration in detail in Section~\ref{sec:firstiter} before supplying the general formulas in Section~\ref{sec:genforms}. 
Algorithm~\ref{alg:IRMS} contains pseudocode for the illustrated  procedure. 

\subsection{First Iteration} 
\label{sec:firstiter}

For motivation, we illustrate the first iteration of Algorithm~\ref{alg:IRMS} in detail. 

Fix a block size $b$, and let $\mathbf{\Omega}^{(0)} \in \F^{n \times b}$ be a random matrix. 
Form the random sample $\mY^{(0)} = \ma \mathbf{\Omega}^{(0)} \in \F^{m \times b}$ and compute
\begin{align}
    \label{eq:LU1}
    \begin{array}{*{9}{c}}
           &      &   &  &   \\
       \mP^{(0)}  & \mY^{(0)} & = & \mL^{(0)} & \mU_1^{(0)}. \\
       m \times m & m \times b &  & m \times b & b \times b 
    \end{array}
\end{align}
Define $\mL_1^{(0)} = \mL^{(0)}(1:b,:)$ and $\mL_2^{(0)} = \mL^{(0)}(b+1:m,:)$.
In practice, the permutation $\mP^{(0)}$ would be stored as a vector (e.g. (\ref{eq:QR})); however, we present the formulas throughout this section in terms of permutation matrices for demonstration purposes.
Let $I_s = (\mP^{(0)} \bmat{1 & ... & m}^*)(1:b)$ be the selected skeleton indices.
As in Section~\ref{sec:errest}, we aim to estimate the approximation error with these skeletons, using an independent random sample of $\ma$.

Let $\mathbf{\Omega}^{(1)} \in \F^{n \times b}$ be an independent random matrix drawn from the same distribution as $\mathbf{\Omega}^{(0)}$, and form $\mY^{(1)} = \ma \mathbf{\Omega}^{(1)} \in \F^{m \times b}$.
Our goal is to estimate 
\begin{align}
    \| \ma - \mW^{(0)} \mR^{(0)} \|_F,
\end{align}
where $\mR^{(0)} = \ma(I_s^{(0)},:)$ and $\mW^{(0)} = (\mP^{(0)})^* \bmat{\mI_{b \times b} \\ \mL_2^{(0)} (\mL_1^{(0)})^{-1}}$.
In computing the error for the rank-$b$ approximation in this iteration, we then have
\begin{align}
\begin{split}
    \label{eq:Aschur}
    \|\ma - \mW^{(0)} \mR^{(0)}\|_F &= \|\mP^{(0)} \ma - \mP^{(0)} \mW^{(0)} \mR^{(0)}\|_F  \\
    &= \left \| \mP^{(0)} \ma - \bmat{(\mP^{(0)} \ma)(1:b,:) \\ \mL^{(0)}_2 (\mL^{(0)}_1)^{-1} (\mP^{(0)} \ma)(1:b,:)} \right \|_F \\
    &=  \| (\mP^{(0)} \ma)(b+1:end,:)  \\
    &\hspace{5mm}- \mL^{(0)}_2 (\mL^{(0)}_1)^{-1} (\mP^{(0)} \ma)(1:b,:) \|_F.
\end{split}
\end{align}
Thus, using the new independent random sample,
\begin{align}
\label{eq:schurerr1}
\begin{split}
    \|\ma - \mW^{(0)} \mR^{(0)} \|_F &\approx \|\ma \mathbf{\Omega}^{(1)} - \mW^{(0)} \mR^{(0)} \mathbf{\Omega^{(1)}} \|_F \\
    &= \|\mP^{(0)} \mY^{(1)}(b+1:end,:) \\
    &\hspace{5mm} - \mL_2^{(0)} (\mL_1^{(0)})^{-1} (\mP^{(0)} \mY^{(1)})(1:b,:) \|_F \\
    &= E^{(1)}_{Schur},
\end{split}
\end{align}
where we note that $E_{Schur}^{(0)} $ can be initialized as an analogous rank-1 approximation. 
If $E_{Schur}^{(1)} > \tau$, then the algorithm should process the new sample $\mY^{(1)}$ while preserving the leading $b\times b$ LUPP factorization of $\mY^{(0)}$, as follows:
\begin{align}
\label{eq:part2}
    \begin{array}{*{9}{c}}
           &      &   &  &   \\
       \mP^{(1)} & \bmat{\mP^{(0)} \mY^{(0)} & \mP^{(0)} \mY^{(1)}} & = & \mL^{(1)} & \mU^{(1)}, \\
       m \times m & m \times 2b &  & m \times 2b & 2b \times 2b 
    \end{array}
\end{align}
where $\mP^{(1)}$ is a permutation matrix to be determined, and
\begin{align}
\label{eq:partL1U1}
\begin{split}
\mL^{(1)}\mU^{(1)} =
\begin{array}{c}
\\ 
     b \\
      \\ 
     \\ \\ 
     m - b \\
     \\ 
     \\
\end{array} 
\begin{tikzpicture}[baseline=0ex]
\matrix (m)[matrix of math nodes,
nodes in empty cells,
 left delimiter={[},
 right delimiter={]},
 inner sep=1.1pt, 
 column 2/.style={black},
 ampersand replacement=\&] 
{
  \mL_{1}^{(0)}     \&    \\     
 \wh{\mL}_{2}^{(0)}   \& \wh{\mL}_{1}^{(1)}   \\  
 \textcolor{green!20}{\wh{\mL}_{2}^{(1)}}   \& \textcolor{red!35}{\wt{\mL}_{22}^{(1)}}    \\ 
  \&      \\  
   \&      \\   
  \textcolor{green!20}{\wt{\mL}_{2}^{(0)}} \&   \wh{\mL}_{2}^{(1)} \\    
     \&      \\  
 \textcolor{green!20}{\mL_{21}^{(1)}}   \& \textcolor{red!35}{\mL_{22}^{(1)}}    \\  
};
\draw[thick, black] (m-1-1.south east)--(m-2-1.south east);
\draw[thick, black] (m-1-1.south east)--(m-1-1.south west);
\draw[thick, black] (m-2-1.south east)--(m-8-1.south east);
\begin{pgfonlayer}{background}
        \draw[draw, fill=red!70, opacity=0.5]
        (m-1-1.south east)--(m-8-1.south east)--(m-8-2.south east)--(m-2-2.south east)--cycle; 
        \draw[draw, fill=gray!40,opacity=0.5]
        (m-1-1.north west)--(m-1-1.south west)--(m-1-1.south east)--cycle;
        \draw[draw, fill=green!40,opacity=0.5]
        (m-1-1.south west)--(m-1-1.south east)--(m-8-1.south east)--(m-8-1.south west)--cycle;
\end{pgfonlayer} 
\end{tikzpicture} 
\begin{tikzpicture}[baseline=-6.5ex]
\matrix (m) [matrix of math nodes,
nodes in empty cells,
 left delimiter={[},
 right delimiter={]},
 inner sep=1.1pt, 
 column 2/.style={black},
 ampersand replacement=\&] 
{
\mU_{1}^{(0)}     \&   \mU_{2}^{(0)}  \\    
   \& \wh{\mU}_{1}^{(1)}   \\  
};
\draw[thick, black](m-1-2.north west)--(m-1-2.south west);
\draw[thick, black](m-1-2.south west)--(m-1-2.south east);
\begin{pgfonlayer}{background}
        \draw[draw, fill=gray!40,opacity=0.5]
        (m-1-1.north west)--(m-1-1.north east)--(m-1-1.south east)--cycle;
        \draw[draw, fill=blue!50,opacity=0.5]
        (m-1-2.north west)--(m-1-2.south west)--(m-1-2.south east)--(m-1-2.north east)--cycle;
        \draw[draw, fill=red!70,opacity=0.5]
        (m-1-2.south west)--(m-1-2.south east)--(m-2-2.south east)--cycle; 
\end{pgfonlayer}
\end{tikzpicture} 
\end{split}
\end{align} 
where $\wh{\mL}_2^{(0)} $ is equal to $\mL_2^{(0)} \in \F^{(m-b)\times b}$ up to row permutation by $\mP^{(1)}$. Then
(\ref{eq:partL1U1}) is equal to a commensurately partitioned version of the left-hand side of (\ref{eq:part2}), denoted
\begin{align}
\label{eq:PY2}
    \underbrace{\bmat{\mbox{\colorbox{gray!20}{$\mI_{b \times b}$}}  &  \\ & \mbox{\colorbox{red!35}{{$\wh{\mP}^{(1)}$}}}}}_{\mP^{(1)}} 
   \begin{tikzpicture}[baseline=-0.5ex]
\matrix (m)[matrix of math nodes,
nodes in empty cells,
 left delimiter={[},
 right delimiter={]},
 inner sep=1pt, 
 column 2/.style={black},
 ampersand replacement=\&] 
{
   \mP^{(0)}\mY^{(0)}  \&  \mP^{(0)}\mY^{(1)}  \\ 
\&   \\
  \textcolor{green!20}{\mP^{(0)}\mY^{(0)}}    \& \textcolor{red!35}{\mP^{(0)}\mY^{(1)}}    \\
 \textcolor{green!20}{\mP^{(0)}\mY^{(0)}}   \& \textcolor{red!35}{\mP^{(0)}\mY^{1)}}    \\ 
      \mP^{(0)}\mY^{(0)}  \& \mP^{(0)}\mY^{(1)} \\  
 \textcolor{green!20}{\mP^{(0)}\mY^{(0)}}   \& \textcolor{red!35}{\mP^{(0)}\mY^{(1)}}    \\ 
 \textcolor{green!20}{\mP^{(0)}\mY^{(0)}}   \& \textcolor{red!35}{\mP^{(0)}\mY^{(1)}}    \\ 
  \textcolor{green!20}{\mP^{(0)}\mY^{(0)}}    \& \textcolor{red!35}{\mP^{(0)}\mY^{(1)}}    \\
};
\draw[thick, black] (m-1-1.north east)--(m-8-1.south east);
\draw[thick, black] (m-3-1.north west)--(m-3-2.north east);
\begin{pgfonlayer}{background}
        \draw[draw, fill=gray!40, opacity=0.5]
        (m-1-1.north west)--(m-3-1.north west)--(m-3-1.north east)--(m-1-1.north east)--cycle; 
        \draw[draw, fill=green!40, opacity=0.5]
        (m-3-1.north west)--(m-8-1.south west)--(m-8-1.south east)--(m-3-1.north east)--cycle;
        \draw[draw, fill=blue!50, opacity=0.5]
        (m-1-2.north west)--(m-3-2.north west)--(m-3-2.north east)--(m-1-2.north east)--cycle; 
        \draw[draw, fill=red!70, opacity=0.5]
        (m-3-2.north west)--(m-8-2.south west)--(m-8-2.south east)--(m-3-2.north east)--cycle; 
\end{pgfonlayer} 
\end{tikzpicture}.
 \begin{array}{c} b \\  \\ \\ \\ m-b \\  \\ \end{array}
\end{align}
Multiplying out (\ref{eq:partL1U1}) and (\ref{eq:PY2}) and matching terms, we have
\begin{align*}
\begin{split}
    \mbox{\colorbox{green!20}{{$\wh{\mL}_2^{(0)}$}}} &= \wh{\mP}^{(1)} \left (  (\mP^{(0)} \mY^{(0)})(b+1:end,:) \right )= \wh{\mP}^{(1)} \mL_2^{(0)}  \\
    \mbox{\colorbox{blue!20}{$\mU_{2}^{(0)}$}} &= \left (\mL_1^{(0)} \right )^{-1} \left ( (\mP^{(0)}\mY^{(1)})(1:b,:) \right ),
\end{split}
\end{align*}
where we obtain $\wh{\mL}_1^{(1)}, \wh{\mL}_2^{(1)}, \wh{\mU}_1^{(1)},$ and $\wh{\mP}^{(1)}$ from an LUPP factorization of 
\begin{align*}
    \bmat{\mbox{\colorbox{red!35}{$\wh{\mL}^{(1)}$}}, \mbox{\colorbox{red!35}{$\wh{\mU}^{(1)}$}}, \mbox{\colorbox{red!35}{$\wh{\mP}^{(1)}$}}} = \texttt{lu} \underbrace{\left ( (\mP^{(0)}\mY^{(1)})(b+1:end,:) - \mL_2^{(0)} \mU_2^{(0)} \right )}_{\ms^{(1)}},
\end{align*}
with $\wh{\mL}_1^{(1)} = \wh{\mL}^{(1)}(1:b,:)$ and  $\wh{\mL}_2^{(1)} = \wh{\mL}^{(1)}(b+1:end,:)$ as in (\ref{eq:partL1U1}).
We note that $\ms^{(1)}$ is precisely the matrix involved in the error estimate $E_{Schur}^{(1)}$ in (\ref{eq:schurerr1}). 

Once the quantities $\wh{\mL}_1^{(1)}, \wh{\mL}_2^{(1)}, \wh{\mU}^{(1)},$ and $\wh{\mP}^{(1)}$ have been computed via LUPP, 
the matrices are re-blocked for the next iteration according to the following partition:
\begin{align}
\label{eq:part2L1U1}
\begin{split}
\mL^{(1)}\mU^{(1)} = \bmat{\mL_1^{(1)} \\ \mL_2^{(1)}} \mU_1^{(1)} =
\begin{array}{c}
\\ 
     b \\
     b \\ 
     \\ \\ 
     m - 2b \\
     \\ 
     \\
\end{array} 
\begin{tikzpicture}[baseline=0ex]
\matrix (m)[matrix of math nodes,
nodes in empty cells,
 left delimiter={[},
 right delimiter={]},
 inner sep=1.1pt, 
 column 2/.style={black},
 ampersand replacement=\&] 
{
  \mL_{1}^{(0)}     \&    \\     
 \wh{\mL}_{2}^{(0)}   \& \wh{\mL}_{1}^{(1)}   \\  
 \textcolor{white}{\wh{\mL}_{2}^{(1)}}   \& \textcolor{white}{\wt{\mL}_{22}^{(1)}}    \\ 
  \&      \\  
   \&      \\   
  \wh{\mL}_{2}^{(0)} \&   \wh{\mL}_{2}^{(1)} \\    
     \&      \\  
 \textcolor{white}{\mL_{21}^{(1)}}   \& \textcolor{white}{\mL_{22}^{(1)}}    \\  
};
\draw[thick, black] (m-1-1.north west)--(m-2-1.south west);
\draw[thick, black] (m-1-1.north west)--(m-2-2.south east);
\draw[thick, black] (m-2-1.south west)--(m-2-2.south east);
\draw[thick, black] (m-2-1.south west)--(m-8-1.south west);
\draw[thick, black] (m-2-2.south east)--(m-8-2.south east);
\draw[thick, black] (m-8-1.south west)--(m-8-2.south east);
\begin{pgfonlayer}{background}
        \draw[draw, fill=gray!40,opacity=0.5]
        (m-1-1.north west)--(m-2-1.south west)--(m-2-2.south east)--cycle;
        (m-1-1.south west)--(m-1-1.south east)--(m-8-1.south east)--(m-8-1.south west)--cycle;
\end{pgfonlayer} 
\end{tikzpicture} 
\begin{tikzpicture}[baseline=-6.5ex]
\matrix (m) [matrix of math nodes,
nodes in empty cells,
 left delimiter={[},
 right delimiter={]},
 inner sep=1.1pt, 
 column 2/.style={black},
 ampersand replacement=\&] 
{
\mU_{1}^{(0)}     \&   \mU_{2}^{(0)}  \\    
   \& \wh{\mU}_{1}^{(1)}   \\  
};
\draw[thick, black](m-1-1.north west)--(m-1-2.north east);
\draw[thick, black](m-1-1.north west)--(m-2-2.south east);
\draw[thick, black](m-2-2.south east)--(m-1-2.north east);
\begin{pgfonlayer}{background}
        \draw[draw, fill=gray!40,opacity=0.5]
        (m-1-1.north west)--(m-1-2.north east)--(m-2-2.south east)--cycle;
\end{pgfonlayer}
\end{tikzpicture}.
\end{split}
\end{align} 
We also update $\mP^{(1)} = \mP^{(1)} \mP^{(0)}$, and $\mY^{(1)} = \bmat{\mY^{(0)} & \mY^{(1)}}$, so that we have a concise version of (\ref{eq:part2}): $\mP^{(1)} \mY^{(1)} = \mL^{(1)} \mU^{(1)}$.

We now generalize these formulas in the next section.

\subsection{General Formulas}
\label{sec:genforms}

Assume for $t \geq 1$ we have a partially pivoted LU factorization of the form
\begin{align*}
    \label{eq:basegen}
    \begin{array}{ccccccc}
    \mP^{(t-1)} & \mY^{(t-1)} &  = & \mL^{(t-1)} & \mU^{(t-1)} & = & \bmat{\mL_1^{(t-1)} \\ \mL_2^{(t-1)}} \mU_1^{(t-1)} \\
    m \times m & m \times tb &  & m \times tb & tb \times tb
    \end{array}
\end{align*}
To compute $E_{Schur}^{(t)}$, we form the independent random sample $\mY^{(t)} = \ma \mathbf{\Omega}^{(t)} \in \F^{m \times b}$ and the Schur complement
\begin{align*}
    \ms^{(t)} = (\mP^{(t-1)}\mY^{(t)})(tb+1:end,:) - \mL^{(t-1)}_2 \mU^{(t-1)}_2,
\end{align*}
where 
$\mU^{(t-1)}_2 = (\mL_1^{(t-1)})^{-1}(\mP^{(t-1)} \mY^{(t)})(1:tb,:)$. Then $E_{Schur}^{(t)} = \|\ms^{(t)} \|_F$.

It is worth highlighting that $\rbr{E_{Schur}^{(t)}}^2$ serves as an unbiased estimate for the squared Frobenius norm of the corresponding true residual matrices: with $\Wb \in \F^{m \times tb}$ and $\Rb \in \F^{tb \times n}$ after $t$ iterations (\ie, $(t-1)$ loops), $\ms^{(t)} \in \F^{(m-bt) \times b}$ satisfies that
\begin{align*}
    \E_{\Omegab}\sbr{\|\ms^{(t)} \|_F^2} = \nbr{\Ab - \Wb\Rb}_F^2.
\end{align*}

Then assuming that $E_{Schur}^{(t)} > \tau$, we apply LUPP to $\ms^{(t)}$:
\begin{align*}
    [\wh{\mL}^{(t)}, \wh{\mU}^{(t)}, \wh{\mP}^{(t)}] = \texttt{lu} \left (\ms^{(t)} \right ),
\end{align*}
and re-block for the next iteration:

\begin{align*}
    \mL_1^{(t)} &= \bmat{\mL_1^{(t-1)} & \mzero \\  (\wh{\mP}^{(t)}\mL_2^{(t-1)})(1:b,:) & \wh{\mL}^{(t)}(1:b,:)} \\
    \mL_2^{(t)} &=\bmat{\wh{\mP}^{(t)}\mL_2^{(t-1)}(b+1:end,:) & \wh{\mL}^{(t)}(b+1:end,:)}  \\
    \mL^{(t)} &= \bmat{\mL_1^{(t)} \\ \mL_2^{(t)}} \\
    \mU^{(t)} &= \mU_1^{(t)} = \bmat{\mU_1^{(t-1)} & \mU_2^{(t-1)} \\ \mzero & \wh{\mU}^{(t)}} \\
    \mP^{(t)} &= \bmat{\mI_{tb \times tb} & \\ & \wh{\mP}^{(t)}} \mP^{(t-1)} \\
    \mY^{(t)} &= \bmat{\mY^{(t-1)} & \mY^{(t)}}.
\end{align*}

The pseudocode for this procedure is given in Algorithm~\ref{alg:IRMS}, where each iteration (loop) takes $O\rbr{b\rbr{\nnz\rbr{\Ab} + m(bt) + (bt)^2}}$ operations asymptotically. With $t$ iterations (\ie, $(t-1)$ loops), the total asymptotic complexity of the algorithm is given by
\begin{align*}
    O\rbr{\nnz\rbr{\Ab}(bt) + m(bt)^2 + (bt)^3} = O\rbr{\nnz\rbr{\Ab}\abbr{I_s} + m \abbr{I_s}^2 + \abbr{I_s}^3},
\end{align*}
where $\abbr{I_s} = k_t \coloneqq bt$ is the final rank of the output row ID $\Wb\Rb$ determined adaptively according to the error tolerance $\tau$.

\remark{ \label{sec:remark2}
The amount of ``oversampling'' in each iteration of the adaptive algorithm is equal to the block size $b$, where $20 \leq b \leq 100$ is usually sufficient. 
As we demonstrate in our numerical experiments, larger block sizes are better for optimally leveraging GPU-capabilities; however, this may result in more computations than necessary to achieve a certain error tolerance. 

The error estimates involving $\mU_r \in \real^{p \times p}$ only require oversampling by $p$ in each iteration, and the random sample $\mY_r$ required for the computation of $\mU_r$ may be formed a priori, unlike the random samples $\mY^{(t)}$ formed within each iteration of the algorithm that are ultimately used for the resulting LU factorization.
More specifically, before the first iteration, we form the random sample $\mY_r = \ma \mOmega_r$. 
Within each iteration $t$ of the algorithm, given an LU factorization of the independent random sample $\mY^{(t-1)} \in \real^{m \times tb}$, to compute the error estimate  $E_{Schur}^{(t)}$, we form an $m \times b$ independent random sample $\mY^{(t)}$.
Computing the error estimates in the $t$-th iteration based on $\mU_r \in \real^{p \times p}$ requires a small additional LU factorization of $\bmat{\mY^{(t-1)} & \mY_r}$ in each iteration.
}

\begin{algorithm}
\caption{randLUPPadap}
\label{alg:IRMS}
\begin{algorithmic}[1]
\REQUIRE $\ma \in \F^{m \times n}$, block size $b = O(1)$, error tolerance $\tau$
\ENSURE Row skeletons $I_s \subset [1,\ldots,m]$ and interpolation matrix $\mW \in \F^{m \times \abbr{I_s}}$ such that $\nbr{\ma - \mW \mR}_F \lesssim \tau$ for $\mR = \ma(I_s,:)$
\STATE Draw independent random matrix $\mathbf{\Omega}^{(0)} \in \F^{n \times b}$
\STATE Form random sample $\mY^{(0)} = \ma \mathbf{\Omega}^{(0)} \in \F^{m \times b}$
\STATE Compute $[\mL^{(0)}, \mU_1^{(0)}, \mP^{(0)}] = \texttt{lu}(\mY)$
\STATE Initialize  $\mL_1^{(0)} = \mL^{(0)}(1:b,:)$ and $\mL_2^{(0)} = \mL(b+1:end,:)$
\STATE Initialize $E^{(0)} = \tau + 1$
\STATE $t=1$
\WHILE{$E^{(t-1)} > \tau$ }
\STATE Draw independent random matrix $\mathbf{\Omega}^{(t)} \in \F^{n \times b}$
\STATE $\mY^{(t)} = \ma \mathbf{\Omega}^{(t)} \in \F^{m \times b}$ 
\algocom{$O(\nnz(\Ab)b)$}
\STATE $\mU_2^{(t-1)} = \left (\mL^{(t-1)}_1 \right )^{-1} {\left (\mP^{(t-1)} \mY^{(t)} \right )(1:tb,:)}$
\algocom{$O\rbr{(tb)^2 b}$}
\STATE $\ms^{(t)} = (\mP^{(t-1)}\mY^{(t)})(tb+1:end,:) - \mL^{(t-1)}_2 \mU^{(t-1)}_2$
\algocom{$O\rbr{(m-tb)tb^2}$}
\STATE $E^{(t)} = \|\ms^{(t)}\|_F$
\algocom{$O\rbr{(m-tb)b}$}
\STATE $[\wh{\mL}^{(t)},\wh{\mU}^{(t)}, \wh{\mP}^{(t)}] = \texttt{lu}(\underset{\textcolor{blue}{(m-tb) \times b }}{\ms^{(t)}})$ 
\algocom{$O\rbr{(m-tb)b^2}$}
\STATE $\mL^{(t)}_1 = \underset{\textcolor{blue}{(t+1)b \times (t+1)b}}{\bmat{ \mL_1^{(t-1)} & \mathbf{0}_{tb \times b} \\ (\wh{\mP}^{(t)}\mL_2^{(t-1)})(1:b,:) & \wh{\mL}^{(t)}(1:b,:)}} $
\STATE $\mL^{(t)}_2 = \underset{\textcolor{blue}{(m - (t+1)b) \times (t+1)b}}{\bmat{(\wh{\mP}^{(t)}\mL_2^{(t-1)})(b+1:end,:) & \wh{\mL}^{(t)}(b+1:end,:)}}$
\STATE $\mU_1^{(t)} = \underset{\textcolor{blue}{(t+1)b \times (t+1)b}}{\bmat{\mU_1^{(t-1)} & \mU_2^{(t-1)} \\ \mathbf{0}_{b \times tb} & \wh{\mU}^{(t)}}}$
\STATE $\mP^{(t)} = \underset{\textcolor{blue}{m \times m}}{\bmat{\mI_{tb \times tb} & \mzero \\ \mzero & \wh{\mP}^{(t)}}} \mP^{(t-1)}$
\STATE $t = t+1$
\ENDWHILE
\STATE Set $I_s = (\mP^{(t-1)} \bmat{1 & \ldots & m}^T)(1:(t-1)b)$ 
\STATE Set $\mW = \mP^{(t-1)*}\bmat{\mI_{(t-1)b \times (t-1)b} \\ \mL^{(t-1)}_{2} (\mL^{(t-1)}_{1})^{-1}}$ 
\end{algorithmic}
\end{algorithm}

\section{Numerical Experiments}
\label{sec:numexp}

We present the results of numerical experiments in which we investigated both the accuracy and the computational efficiency of Algorithm~\ref{alg:IRMS}. 

\subsection{Accuracy Comparisons}
\label{sec:acccomp}

In these experiments, we compare the accuracy of Algorithm~\ref{alg:IRMS} to several reference methods in constructing a low-rank approximation to different test matrices $\ma \in \real^{m \times n}$. 
We demonstrate the accuracy of the algorithm for the row ID using several different input matrices and random matrix distributions for the randomized sketching.

For the input matrices, we considered the following synthetic matrices for their spectral behavior, as well as two empirical datasets:
\begin{itemize}
\item \textbf{Fast decay:} $\ma = \mU \mtx{D} \mV^*$, where $\mtx{U} \in \real^{m \times n}$ and $\mtx{V} \in \real^{n \times n}$ have orthonormal columns with Haar measure, via \texttt{qr} on two Gaussian random matrices. The matrix  $\mtx{D}$ is diagonal with entries $d_{ii} = \beta^{(i-1)/(n-1)}$ for $\beta = 10^{-16}$ and $m,n = 5000$.
\item \textbf{Kahan:} \cite{kahan1966} $\ma=\mtx{D}_n\mtx{K}_n\in\real^{n\times n}$, where 
    \begin{align*}
    \mtx{D}_n &= \diag\begin{bmatrix}1 & \zeta & \zeta^2 & \cdots & \zeta^{n-1}
    \end{bmatrix}, \\
        \mtx{K}_n &=\begin{bmatrix} 
        1 & -\varphi & -\varphi & \cdots & -\varphi \\
         & 1 &  -\varphi & \cdots & -\varphi  \\
         &  & \ddots & \ddots & \vdots \\
         &  &  & 1 &  -\varphi  \\
         &  &  &  & 1
        \end{bmatrix},
    \end{align*}
with $\zeta^2 + \varphi^2 = 1$ for $\zeta = 0.99$ and $n = 5000$.
\item \textbf{MNIST}: \cite{deng2012mnist} training set of size $60,000 \times 784$, comprising flattened and normalized images of hand-written digits 0-9.
\item \textbf{large}: \cite{large} The SparseSuite matrix `\texttt{large},' a sparse full-rank $4,282 \times 8,617$ matrix with $20,635$ nonzero entries generated by a linear programming problem sequence.
\end{itemize}

For randomized sketching, we used random matrices from the following distributions:
\begin{itemize}
    \item \textbf{Gaussian}: $\mOmega^{(t)} \in \real^{n \times b}, \mOmega_r \in \real^{n \times p}$ have \textit{i.i.d} Gaussian entries with standard deviations $1/b$ and $1/p$, resp
    \item \textbf{Subsampled Randomized Trigonometric Transforms (SRTT) }\cite{halko2011, tropp2011}: $\mOmega^{(t)}, \mOmega_r = \sqrt{\frac{n}{\ell}} \mtx{\Pi}_{n} \mtx{T} \mtx{\Phi} \mtx{Pi}_{n \rightarrow \ell}$ for $\ell = b, p$, resp. Here, $\mtx{\Pi}_{n \rightarrow \ell} \in \real^{n \times \ell}$ is a uniformly random selection of $\ell$ out of $n$ columns, whereas $\mtx{\Pi} \in \real^{n \times n}$ is a random permutation of the $n$ columns.
    The matrix $\mtx{T} \in \mathbb{C}^{n \times n}$ is a discrete Fourier transform, and $\mtx{\Phi} \coloneqq \diag(\phi_1, \ldots, \phi_n)$ has $\textit{i.i.d.}$ Rademacher random variable entries that change signs randomly.
    \item \textbf{Sparse Sign} \cite{nelson2013, woodruff2015}: $\mOmega^{(t)} = \sqrt{\frac{n}{\zeta}} \bmat{\mtx{s}_1,\ldots,\mtx{s}_n}^T$ for some $2 \leq \zeta \leq \ell$ (we choose $\zeta =8$), where each row $\mtx{s}_j \in \real^{\ell}$ for $j = 1,\ldots,n$ is filled with $\zeta$ independent Rademacher random variables at uniformly random coordinates.
\end{itemize}

For each of the input matrices and random matrices above, we computed the following  approximation errors and estimates in the Frobenius norm in every iteration $t$ of Algorithm~\ref{alg:IRMS}:
\begin{itemize}
    \item \textbf{SVD}: $(\sum_{k_t+1}^{r} \sigma_j^2)^{1/2}$ for $r = \min(m,n)$ and approximation rank $k_t \equiv tb$,
    \item \textbf{\texttt{randLUPPadap}-\textbf{Schur}}: (Algorithm~\ref{alg:IRMS})  $E^{(t)}_{Schur} = \|\ms^{(t)}\|_F$ as in Section~\ref{sec:genforms},
    \item \textbf{\texttt{randLUPPadap}-ID}: if $P_{LU}$ denotes the row permutation vector resulting from \texttt{randLUPPadap}, we compute the row ID error
    \begin{align*} 
    \left \| \ma(P_{LU},:) - \mW_{LU}^{(t)} \ma(P_{LU}(1:k_t),:) \right \|_F
    \end{align*}
    where $\mW_{LU}^{(t)}$ denotes the row interpolation matrix corresponding to approximation rank $k_t$,
    \item \textbf{\texttt{randLUPPadap}-sID}:  if $P_{LU}$ denotes the row permutation vector resulting from \texttt{randLUPPadap}, we compute the stable row ID error
    \begin{align*}
        \| \ma^* - \mQ_{LU}^{(t)} (\mQ_{LU}^{(t)})^* \ma^* \|_F
    \end{align*}
    where $\mQ_{LU}^{(t)} = \texttt{qr}\left (\ma(P_{LU}(1:k_t),:)^* \right)$ as in (\ref{eq:sID_row}),
    \item \textbf{\texttt{randCPQR}-ID}: \cite{voronin2017} if $P_{QR}$ denotes the row permutation vector resulting from \texttt{randCPQR} applied to $(\mY^{(t)})^*$, we compute the row ID error
    \begin{align*}
        \left \| \ma(P_{QR},:) - \mW_{QR}^{(t)} \ma(P_{QR}(1:k_t),:) \right \|_F \end{align*}
    where $\mW_{QR}^{(t)}$ is the row interpolation matrix from \texttt{randCPQR} given by $\mW_{QR}^{(t)}(P_{QR},:) = \bmat{\mI_{k_t} & \mR_{11}^{-1} \mR_{12}}^*$ as in (\ref{eq:cpqr}), 
    %
    \item \textbf{\texttt{randCPQR}-sID}: if $P_{QR}$ denotes the row permutation vector resulting from \texttt{randCPQR} applied to $(\mY^{(t)})^*$, we compute the stable row ID error
    \begin{align*}
        \left \| \ma^* - \mQ_{QR}^{(t)} (\mQ_{QR}^{(t)})^* \ma^* \right \|_F \end{align*}
    where $\mQ_{QR}^{(t)} = \texttt{qr}\left (\ma(P_{QR}(1:k_t),:)^* \right)$ as in (\ref{eq:sID_row}).
\end{itemize}

For reference, in relation to Remarks~\ref{sec:remark1} and \ref{sec:remark2}, we also plot the quantities in (\ref{eq:normUest}) and (\ref{eq:maxUest}) involving $\mU_r$, denoted by $\|\mU_r\|_F$ and $\max|\mU_r|$, respectively.

\begin{figure}
\centering
\begin{subfigure}{0.49\textwidth}
\centering
\includegraphics[width=\textwidth]{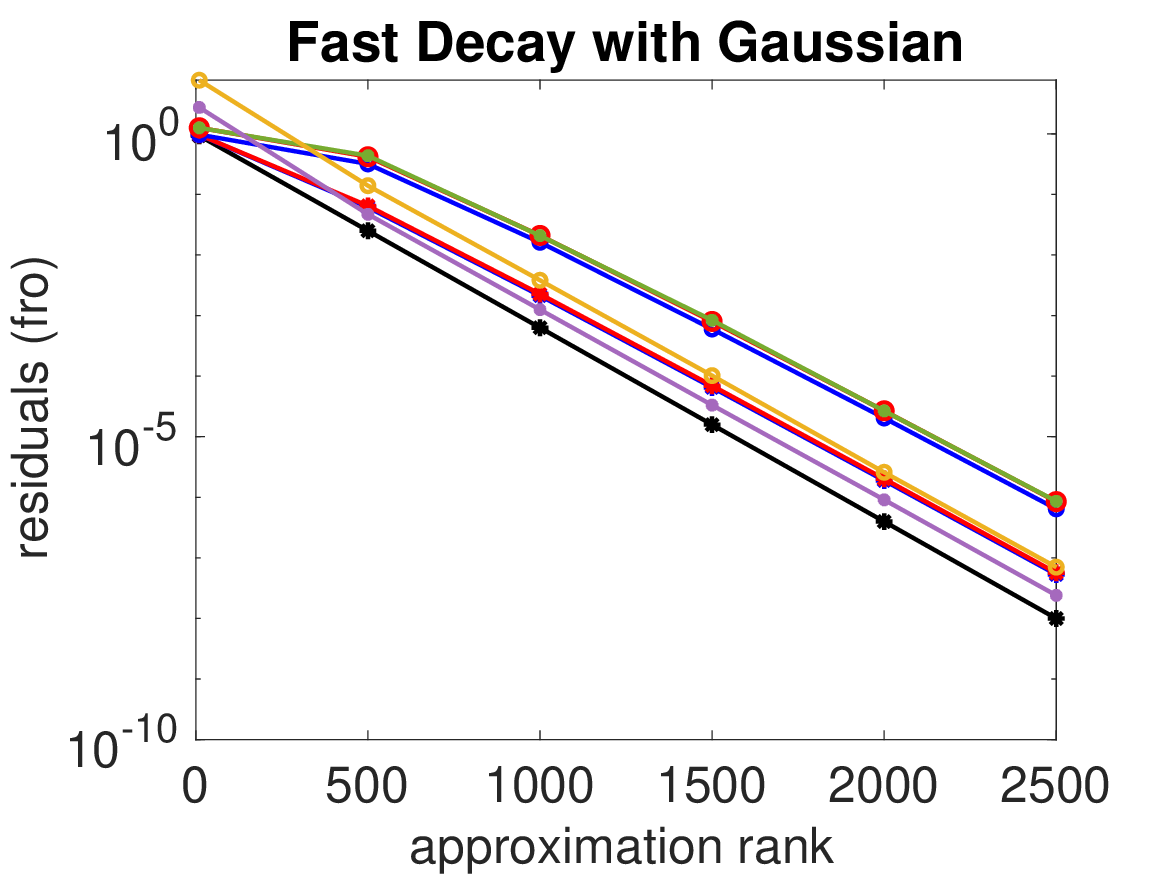}
\caption{}
\end{subfigure} 
\hfill
\begin{subfigure}{0.29\textwidth}
\centering
\includegraphics[width=\textwidth]{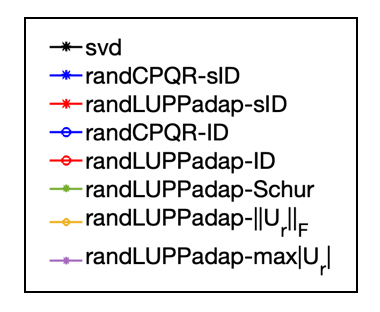}
\vspace{0mm} 
\end{subfigure} 
\begin{minipage}{0.19\textwidth}
    \textcolor{white}{nothing}
\end{minipage}\\
\vspace{3mm}
\begin{subfigure}{0.49\textwidth}
\centering
\includegraphics[width=\textwidth]{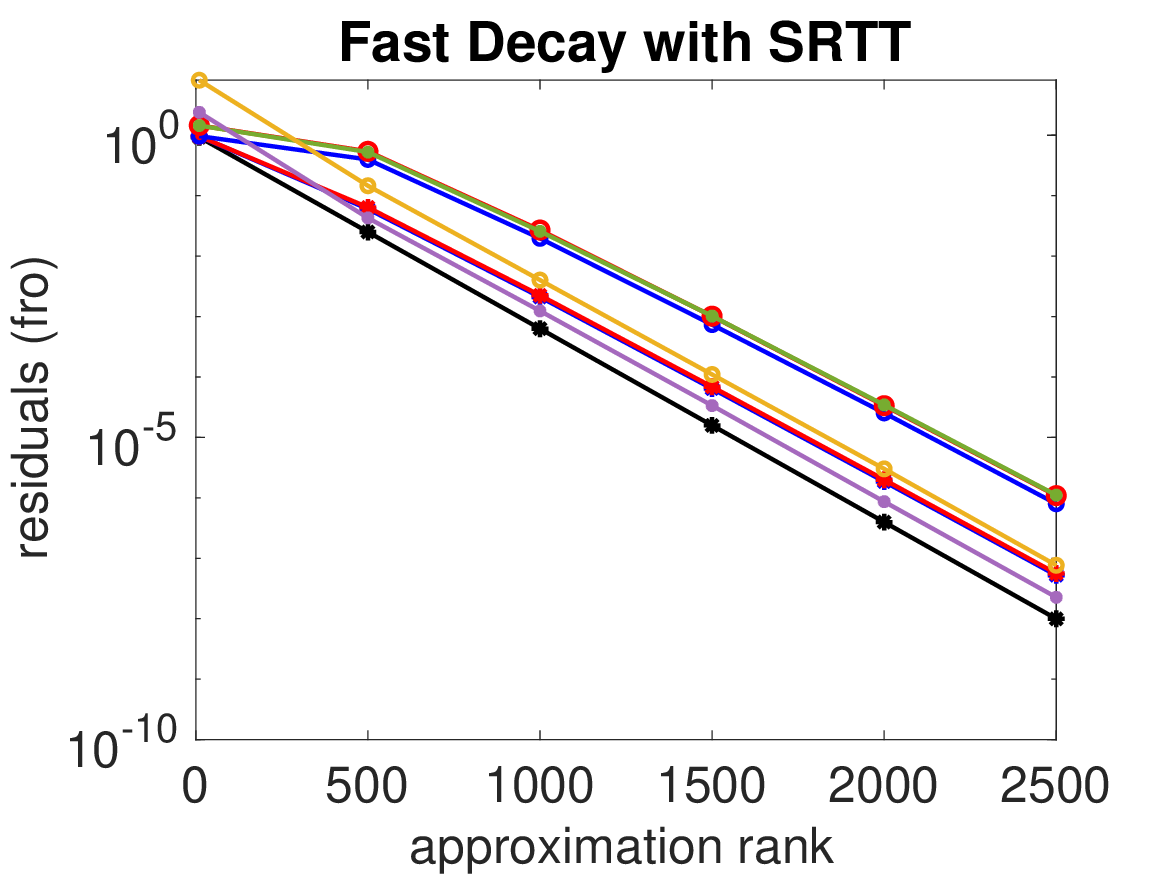}
\caption{}
\end{subfigure}
\begin{subfigure}{0.49\textwidth}
\centering
\includegraphics[width=\textwidth]{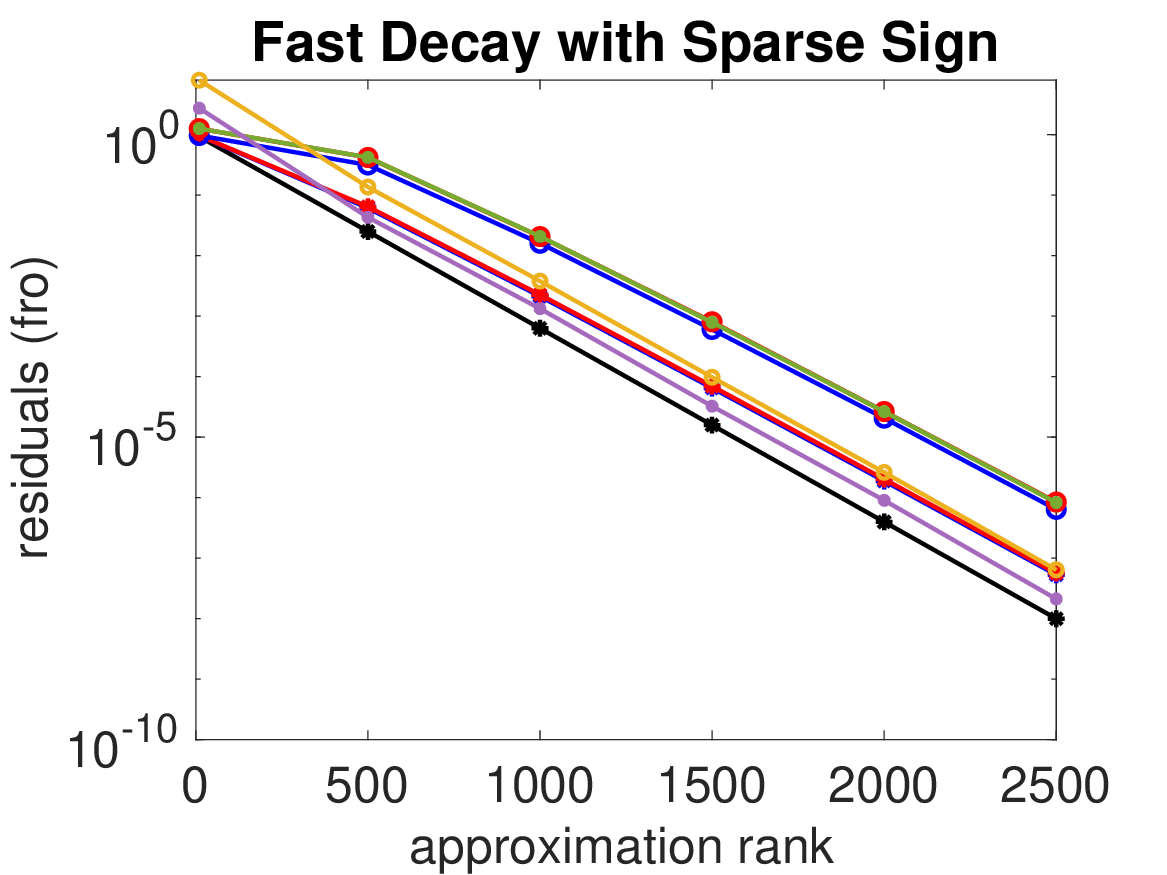}
\caption{}
\end{subfigure}
\caption{Row ID accuracy comparisons for synthetic test matrix \textbf{Fast Decay} in Section~\ref{sec:acccomp}.The mean of 50 realizations of Algorithm~\ref{alg:IRMS} are plotted for each error quantity, with block size $b = 500$ and $p=10$. Results are shown for approximation ranks up to $n/2$, where each marker indicates an iteration of Algorithm~\ref{alg:IRMS}.
}
\label{fig:exp_acc}
\end{figure}

\begin{figure}
\centering
\begin{subfigure}{0.49\textwidth}
\centering
\includegraphics[width=\textwidth]{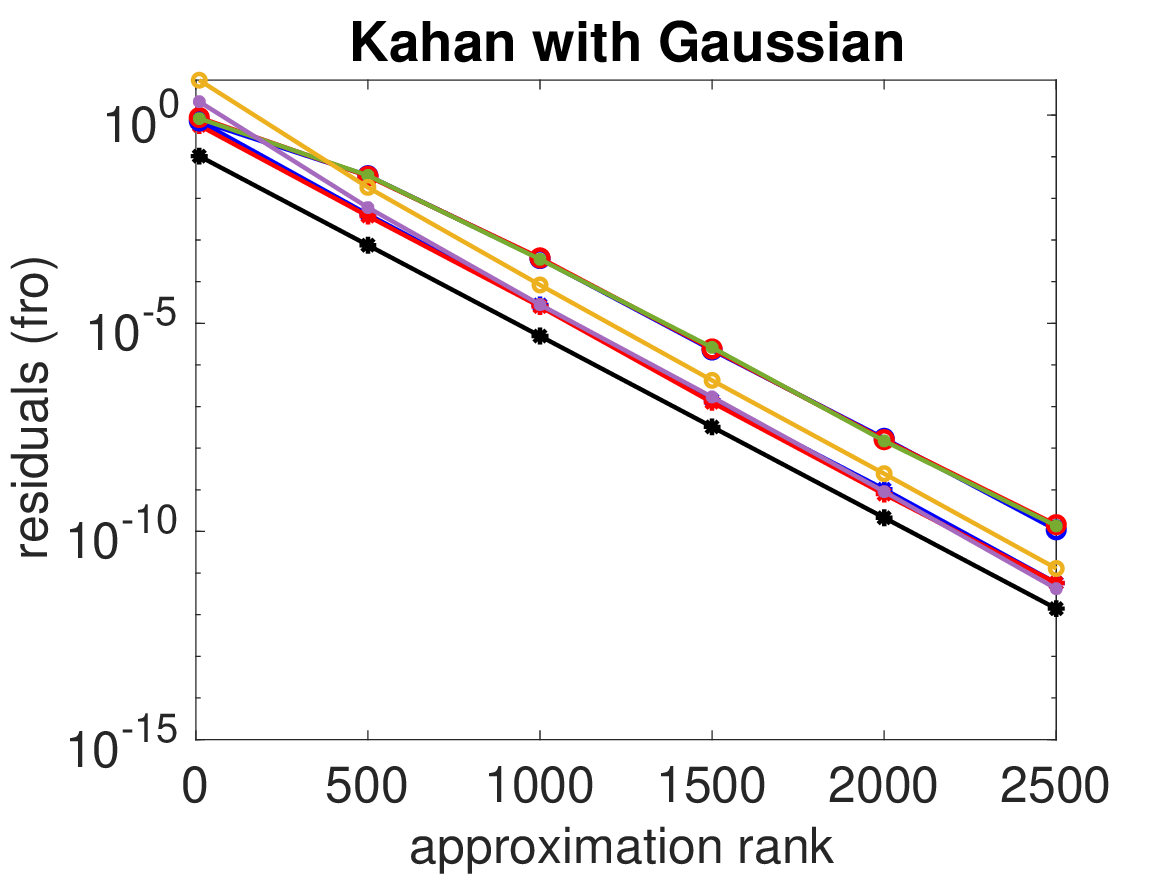}
\caption{}
\end{subfigure} 
\hfill
\begin{subfigure}{0.29\textwidth}
\centering
\includegraphics[width=\textwidth]{AccFigs_FastDecay/legend.png}
\vspace{0mm}
\end{subfigure} 
\begin{minipage}{0.19\textwidth}
    \textcolor{white}{nothing}
\end{minipage}\\
\vspace{3mm}
\begin{subfigure}{0.49\textwidth}
\centering
\includegraphics[width=\textwidth]{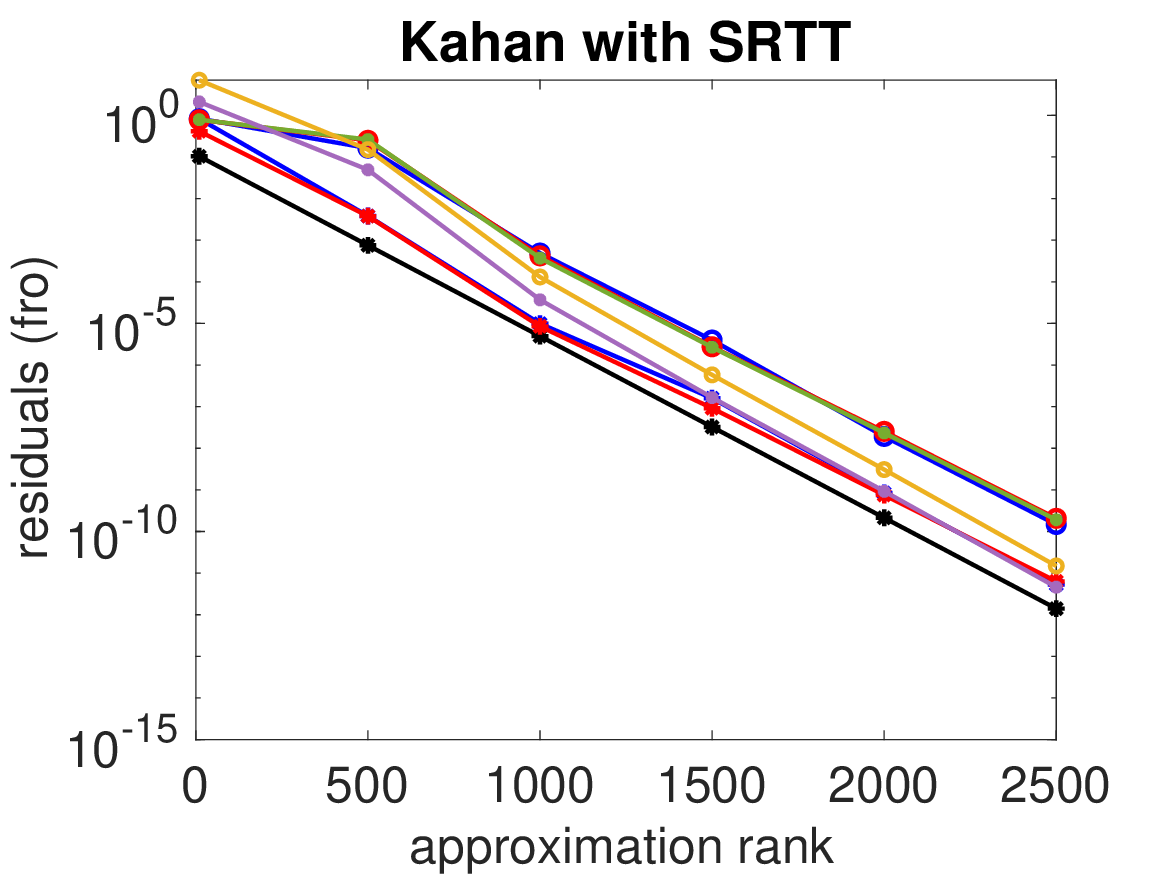}
\caption{}
\end{subfigure}
\begin{subfigure}{0.49\textwidth}
\centering
\includegraphics[width=\textwidth]{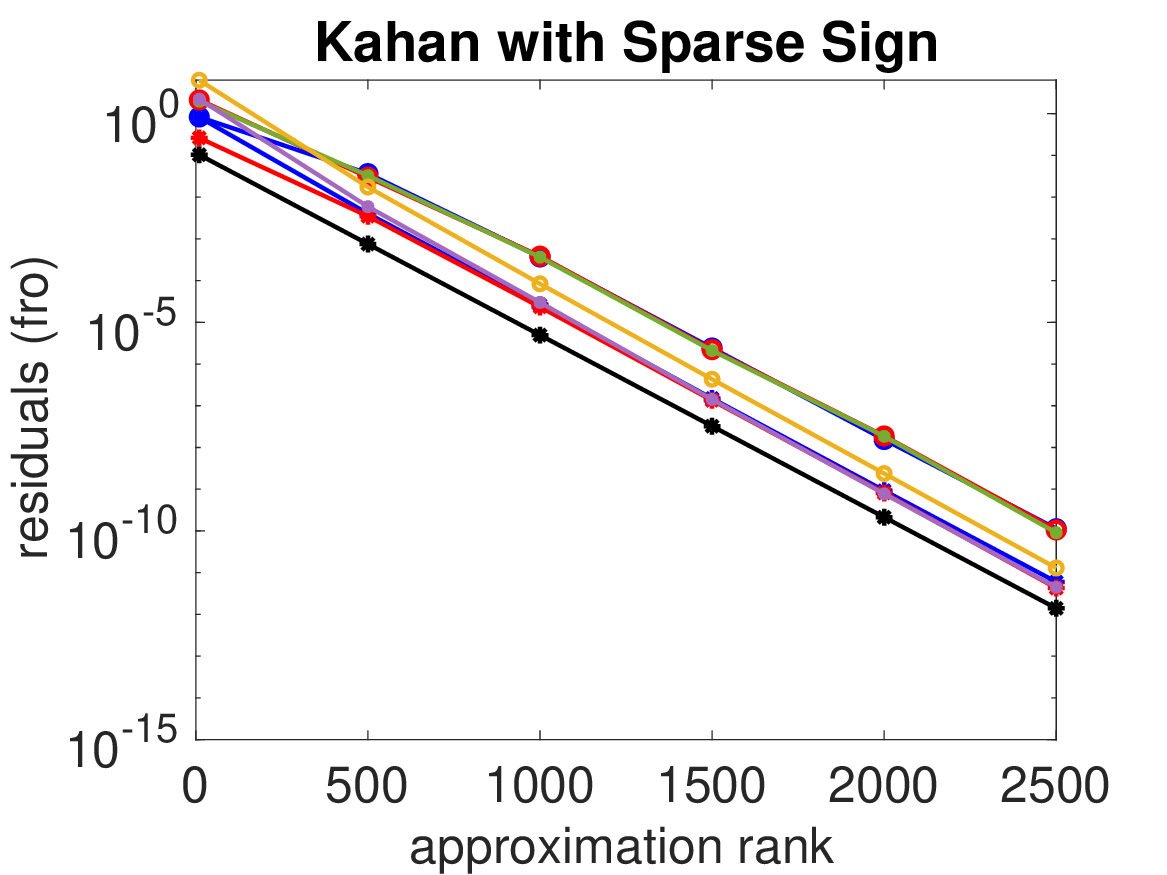}
\caption{}
\end{subfigure}
\caption{Row ID accuracy comparisons for synthetic test matrix \textbf{Kahan} in Section~\ref{sec:acccomp}. The mean of 50 realizations of Algorithm~\ref{alg:IRMS} are plotted for each error quantity, with block size $b = 500$ and $p=10$. Results are shown up to an approximation rank of $n/2$, where each marker indicates an iteration of Algorithm~\ref{alg:IRMS}.
}
\label{fig:kahan_acc}
\end{figure}

\begin{figure}
\centering
\begin{subfigure}{0.49\textwidth}
\centering
\includegraphics[width=\textwidth]{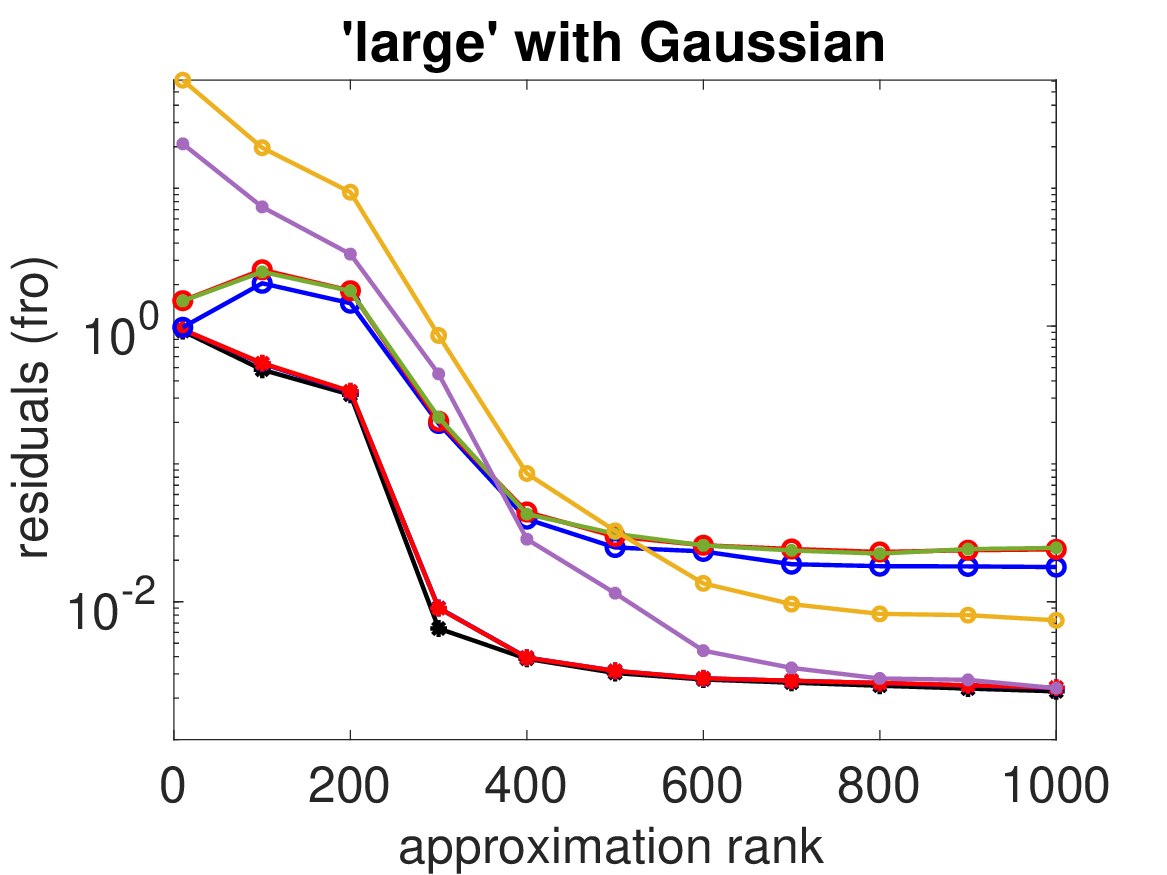}
\caption{}
\end{subfigure} 
\hfill
\begin{subfigure}{0.29\textwidth}
\includegraphics[width=\textwidth]
{AccFigs_FastDecay/legend.png}
\vspace{0mm}
\end{subfigure} 
\begin{minipage}{0.19\textwidth}
    \textcolor{white}{nothing}
\end{minipage}\\
\vspace{3mm}
\begin{subfigure}{0.49\textwidth}
\centering
\includegraphics[width=\textwidth]{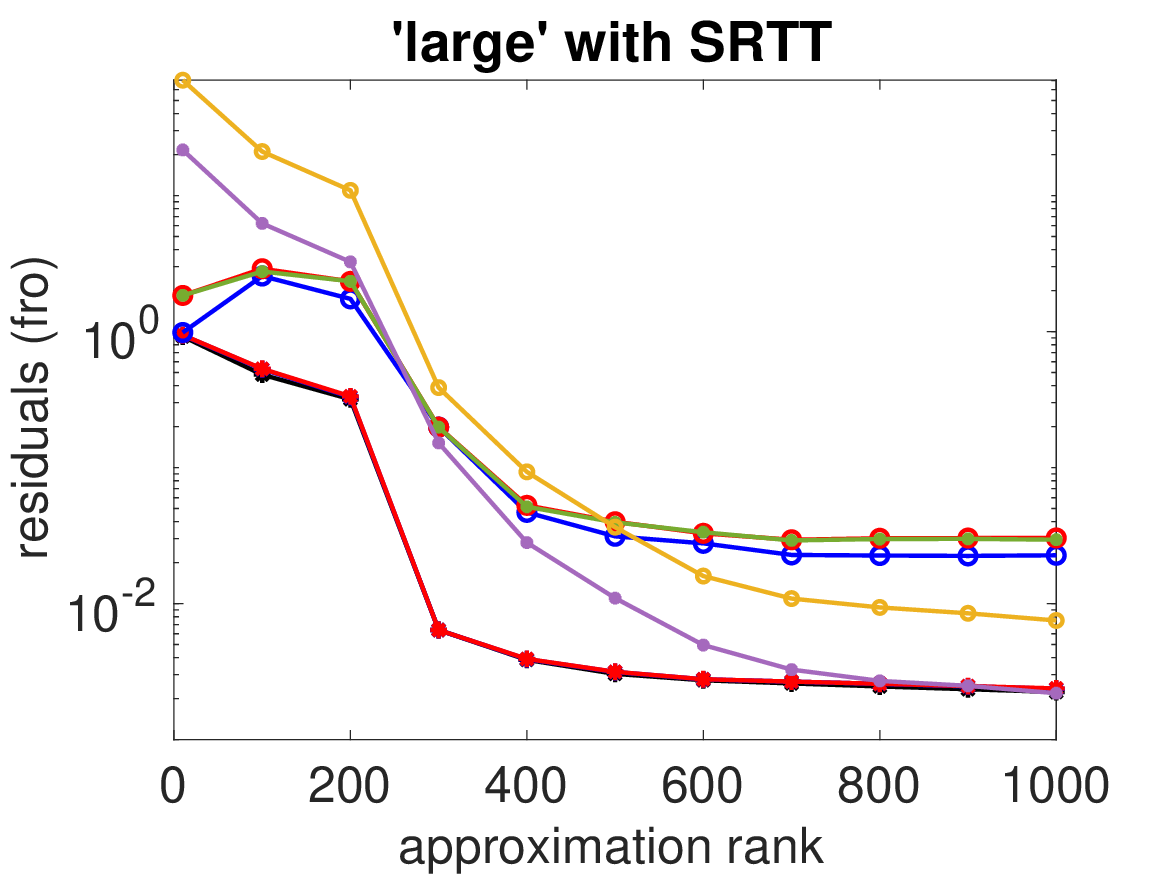}
\caption{}
\end{subfigure}
\begin{subfigure}{0.49\textwidth}
\centering
\includegraphics[width=\textwidth]{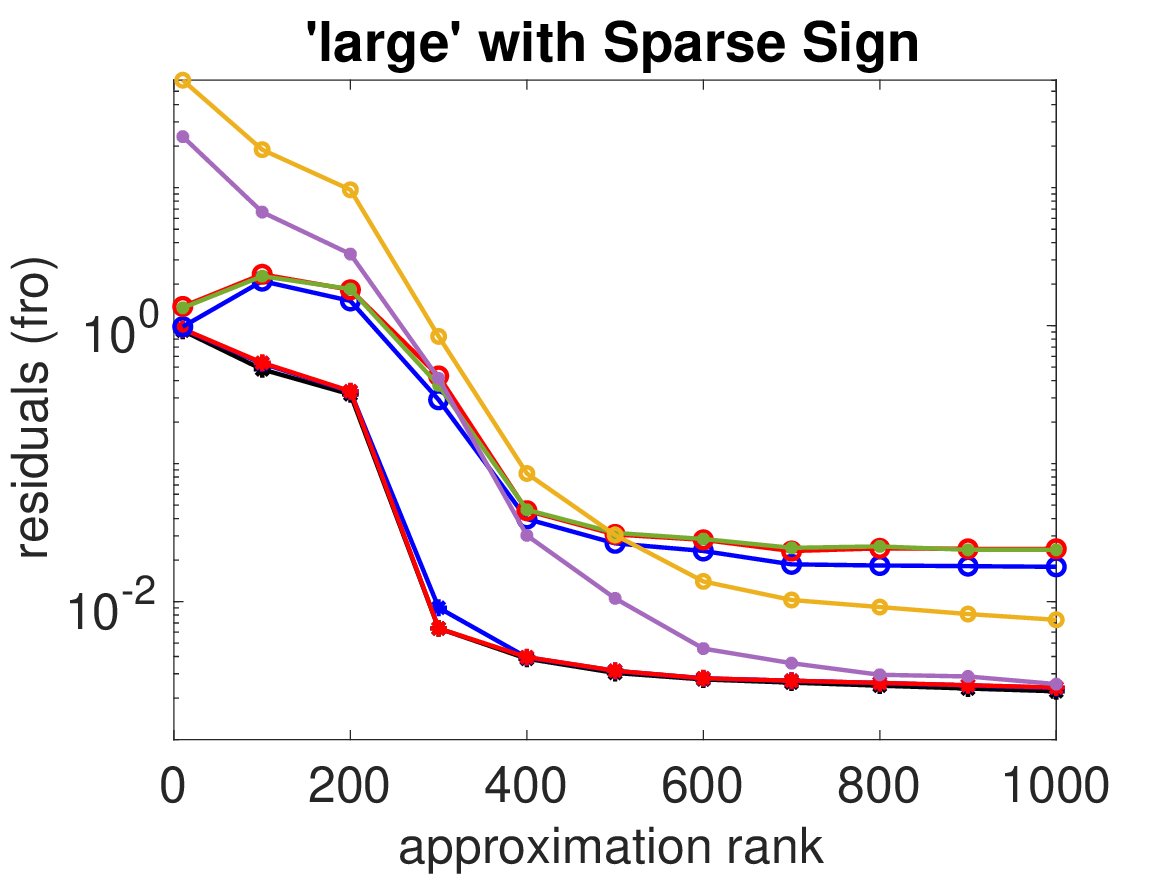}
\caption{}
\end{subfigure}
\caption{Row ID accuracy comparisons for the empirical dataset \textbf{large} in Section~\ref{sec:acccomp}. The mean of 50 realizations of Algorithm~\ref{alg:IRMS} are plotted for each error quantity, with block size $b = 100$ and $p=10$. Each marker indicates an iteration of Algorithm~\ref{alg:IRMS}.
}
\label{fig:large_acc}
\end{figure}

\begin{figure}
\centering
\begin{subfigure}{0.49\textwidth}
\centering
\includegraphics[width=\textwidth]{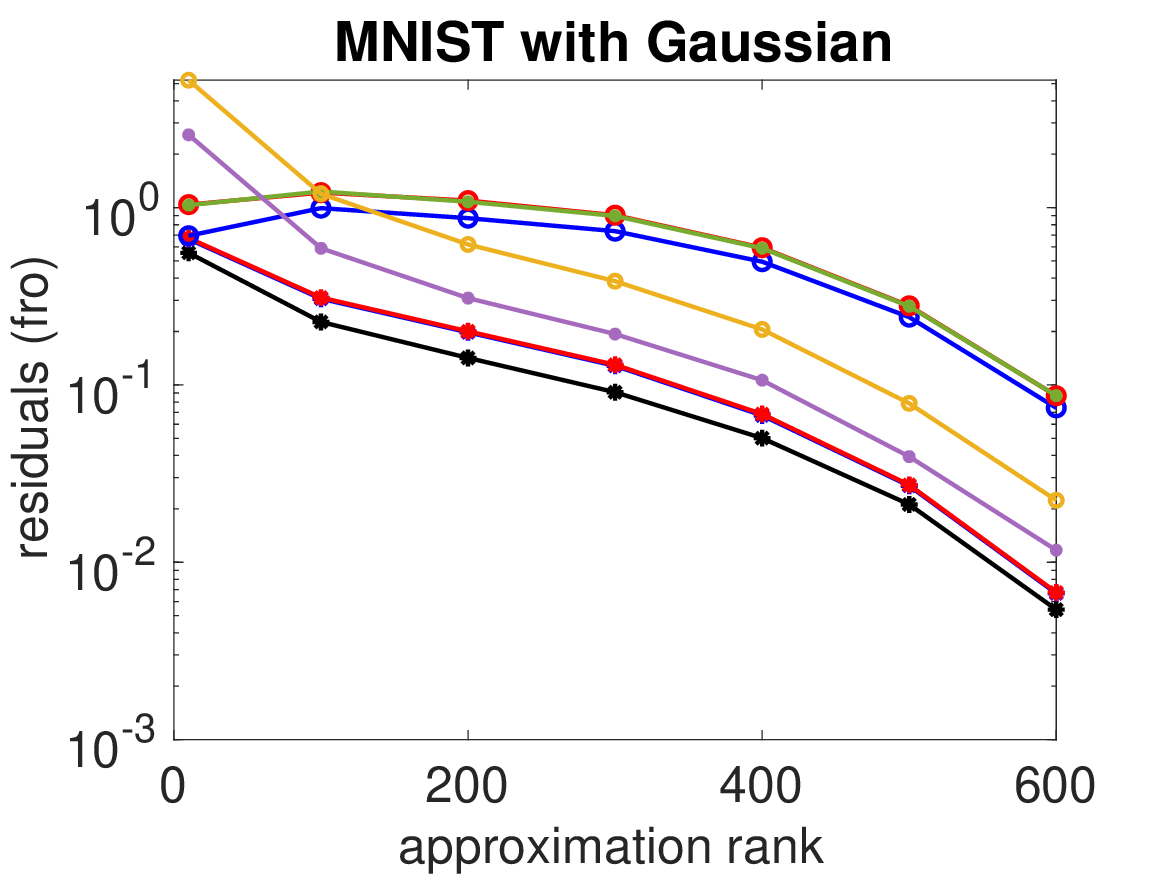}
\caption{}
\end{subfigure} 
\hfill
\begin{subfigure}{0.29\textwidth}
\includegraphics[width=\textwidth]{AccFigs_FastDecay/legend.png}
\vspace{0mm}
\end{subfigure} 
\begin{minipage}{0.19\textwidth}
    \textcolor{white}{nothing}
\end{minipage}\\
\vspace{3mm}
\begin{subfigure}{0.49\textwidth}
\centering
\includegraphics[width=\textwidth]{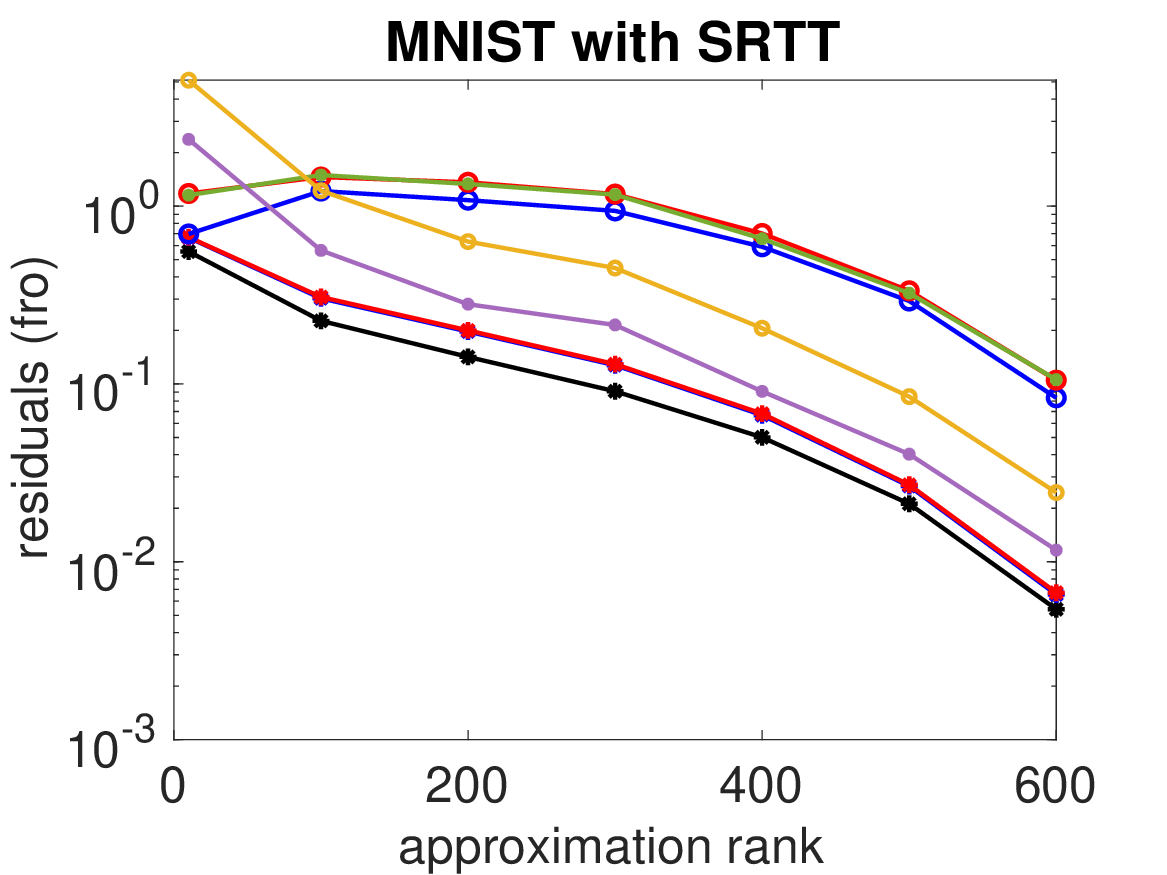}
\caption{}
\end{subfigure}
\begin{subfigure}{0.49\textwidth}
\centering
\includegraphics[width=\textwidth]{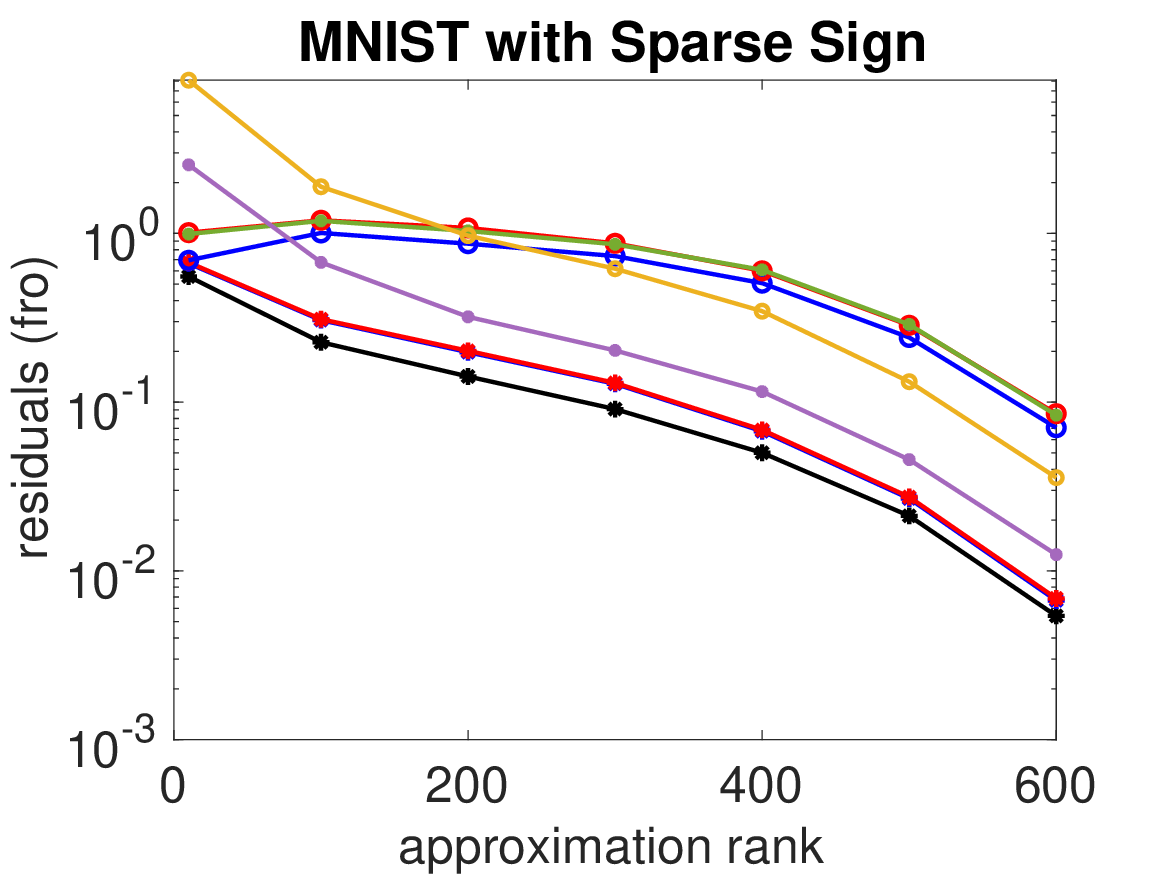}
\caption{}
\end{subfigure}
\caption{Row ID accuracy comparisons for the empirical dataset \textbf{MNIST} in Section~\ref{sec:acccomp}. The mean of 50 realizations of Algorithm~\ref{alg:IRMS} are plotted for each error quantity, with block size $b = 100$ and $p=10$. Each marker indicates an iteration of Algorithm~\ref{alg:IRMS}.
}
\label{fig:mnist_acc}
\end{figure}

In Figures~\ref{fig:exp_acc} and \ref{fig:kahan_acc}, we plot the residual errors and error estimates after each iteration $t$ of \texttt{randLUPPadap} for the synthetic test matrices \textbf{Fast Decay} and \textbf{Kahan}, respectively.
In Figures~\ref{fig:large_acc} and \ref{fig:mnist_acc}, we plot the residual errors for the empirical datasets \textbf{large} and \textbf{MNIST}.

The experiments confirm our claim in Section~\ref{sec:genforms} that $\|\ms^{(t)}\|_F^2 = (E^{(t)}_{Schur})^2$ is an unbiased estimator for the row ID error $\| \ma - \mW \mR \|_F^2$ in \texttt{randLUPPadap}.
This is evidenced in each figure by the overlapping lines labeled {``randLUPPadap-Schur''} and {``randLUPPadap ID''}. 

As expected, both \texttt{randLUPPadap} and \texttt{randCPQR} display comparable accuracy when the errors are computed with the stable ID method as in (\ref{eq:sID_row}), which gives the true ID approximation error with the selected skeletons. 
Because \texttt{randLUPPadap} is mathematically equivalent to \texttt{randLUPP} with the same overall embedding dimension, these results corroborate the findings of \cite{dong2022} that the stable ID computations for randomized LUPP yield the same quality of error as stable ID computations for randomized CPQR.

We also observe comparable accuracy in the other computed error estimates without orthonormalization. 
These close quantities (\texttt{randLUPPadap}-Schur, {\texttt{randLUPPadap}-ID}, {\texttt{randCPQR}-ID}) indicate that the accuracy of Algorithm~\ref{alg:IRMS} in selecting skeletons is competitive with \texttt{randCPQR}.
Additionally, when there is fast spectral decay like in Figures~\ref{fig:exp_acc} and \ref{fig:kahan_acc}, the Schur complement estimate used in Algorithm~\ref{alg:IRMS} will achieve the same degree of accuracy as stable ID errors or singular value lower bounds, if they were used instead for the algorithm's termination criteria. 
For the empirical datasets in Figure~\ref{fig:large_acc} and \ref{fig:mnist_acc} with slower spectral decay, more oversampling would be required for the plain row ID errors to achieve the same level of accuracy as the stable ID errors. 
We note that in all cases, the error estimates involving the residual upper triangular matrix $\mU_r$ are more accurate than the error estimates with more oversampling ($b$ vs. $p$) than in the computation of $\mU_r$.
In general, the estimates with $\|\mU_r\|_F$ and $\max|\mU_r||$ closely approximate the stable ID errors computed for both \texttt{randLUPPadap} and \texttt{randCPQR}, and our findings support the claim in Remark~\ref{sec:remark1} that $\|\mU_r\|_F$ can provide a reliable ID error estimate (i.e. does not yield underestimates as $\max|\mU_r|$ can in rare cases).

In summary:

\begin{enumerate}
\item
The proposed error estimator $\|\ms^{(t)}\|_F$ is an efficient and accurate estimate of the true ID error. 

\item 
The ID errors of all methods are close (so are the \emph{stable} ID errors). This implies that the new method based on LUPP is faster without sacrificing accuracy.

\item
The ID error is usually larger than the \emph{stable} ID error, especially when singular values decay slowly. The \emph{stable} ID error shows the true accuracy of the skeletons.

\end{enumerate}

\subsection{Speed and Performance Results}

In this section, we compare  the  new method \texttt{randLUPPadap} to two reference methods \texttt{randLUPP} and \texttt{randCPQR} in terms of their  running time for computing IDs of the  \textbf{Fast Decay} matrices as described in Section~\ref{sec:acccomp}. 
The machine where our numerical experiments were performed has
\begin{itemize}

\item
an NVIDIA GPU that is a Tesla V100 GPU with 32 GB of memory (peak performance $\approx$ 7 TFlop/s in double precision),

\item
two Intel Xeon Gold 6254 CPUs, each with 18 cores at 3.10 GHz (peak performance $\approx$ 1.27 TFlop/s in double precision),

\item
{a PCIe 3.0 $\times 16$ between the CPU and the GPU that can deliver up to 15.75 GB/s.}

\end{itemize}
Our GPU code was compiled with the NVIDIA compiler {nvcc} (version 11.8.89) and linked with the cuBLAS library (version 11.11.03) on the Linux OS (5.15.0-60-generic.x86\_64). Our CPU code was compiled with the GNU compiler gcc (version 11.3.0) and linked with the multi-threaded Intel MKL library (version 2022.1.0). Our CPU code was ran with 32 threads (OMP\_NUM\_THREADS=32), which empirically led to its best performance. All calculations were performed with double-precision floating-point arithmetic. 
Timing results were taken as the average of five consecutive runs.

The two reference methods are implemented as follows:
\begin{itemize}

\item \texttt{randLUPP}: the GPU version calls the \texttt{getrf} routine from the cuBLAS library, and the CPU version calls the \texttt{getrf} routine from the multi-threaded Intel MKL library.

\item \texttt{randCPQR}: the GPU version calls the \texttt{geqp3\_gpu} routine from the MAGMA library~\cite{dghklty14} (no \texttt{geqp3} in cuBLAS), and the CPU version calls the \texttt{geqp3} routine from the multi-threaded Intel MKL library.

\end{itemize}

\begin{figure}
\centering
\begin{subfigure}{0.49\textwidth}
\centering
\includegraphics[width=\textwidth]{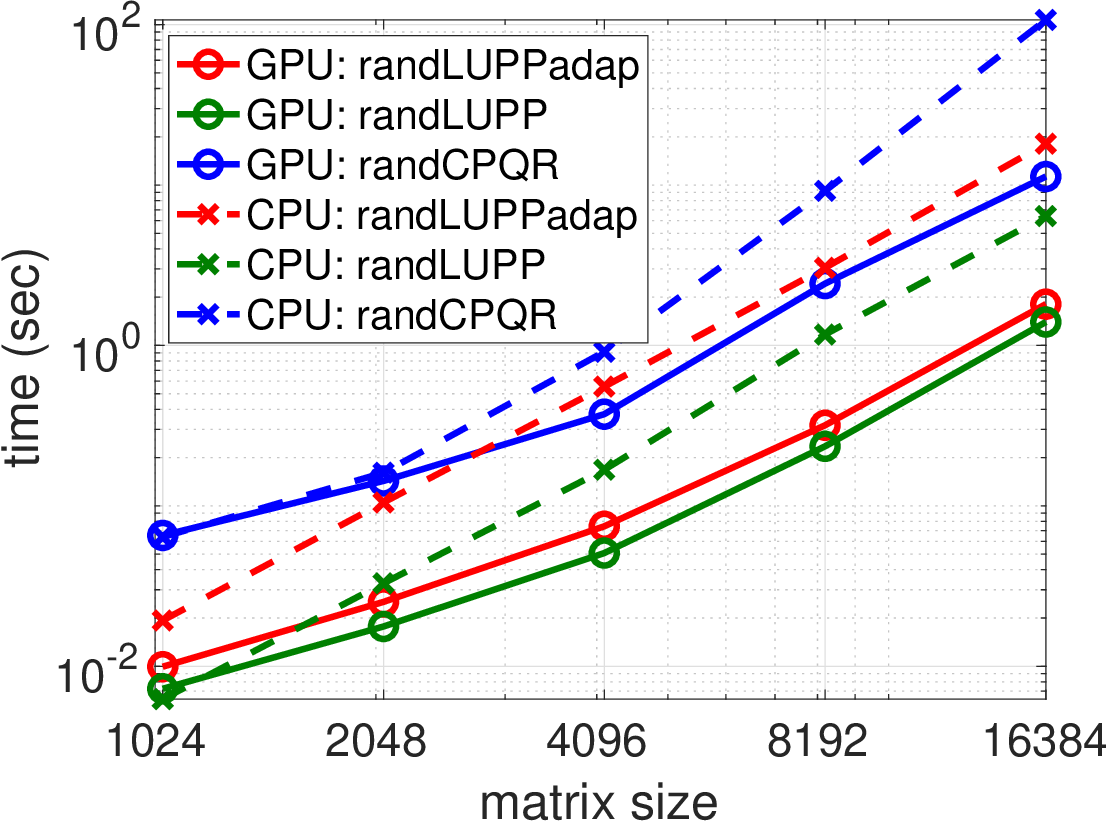}
\caption{Tolerance $\tau = 10^{-8}$}
\end{subfigure}
\begin{subfigure}{0.49\textwidth}
\centering
\includegraphics[width=\textwidth]{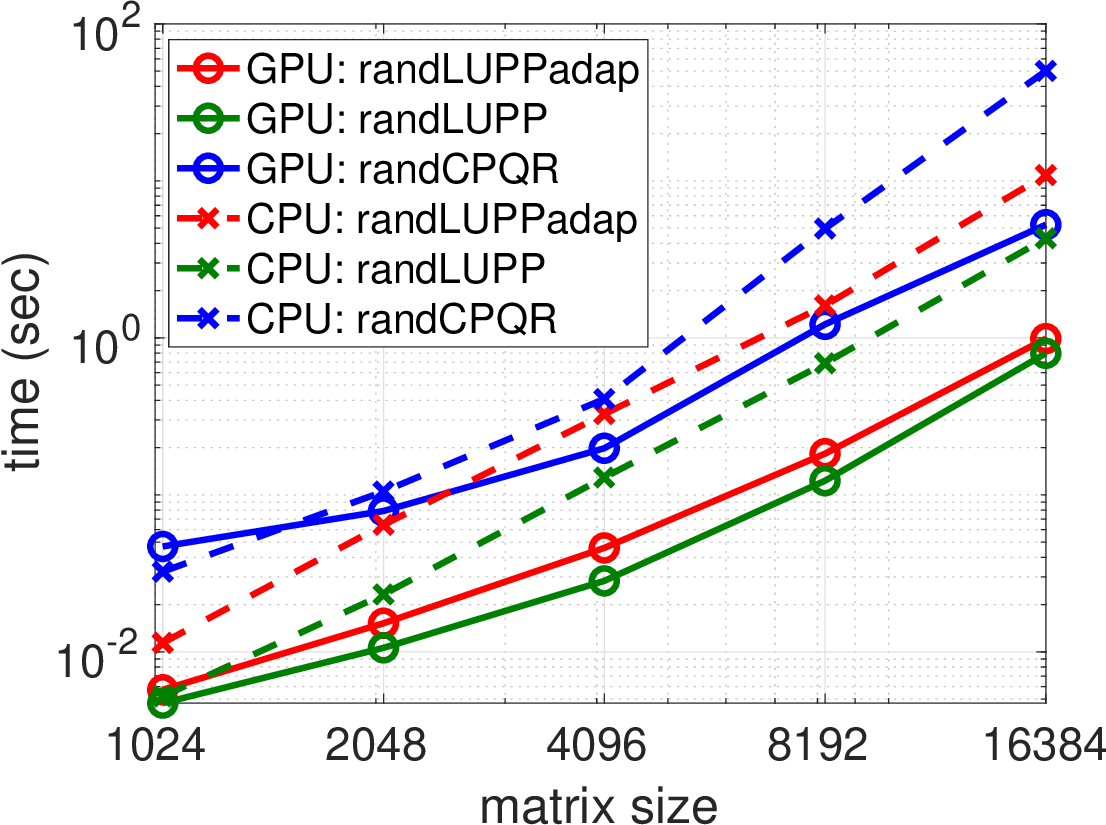}
\caption{Tolerance $\tau = 10^{-4}$}
\end{subfigure}
\caption{Timings for computing the IDs of $n \times n$ \textbf{Fast Decay} matrices as described in Section~\ref{sec:acccomp}.  The two tolerances $10^{-8}$ and $10^{-4}$ correspond to a numerical rank of approximately $n/2$ and $n/4$, respectively. The block size $b$ used in \texttt{randLUPPadap}  is 128 by default. The computed numerical rank from \texttt{randLUPPadap}  is used as an input for the two reference methods \texttt{randLUPP} and \texttt{randCPQR}. 
}
\label{f:time}
\end{figure}

\begin{figure}
\centering
\begin{subfigure}{0.49\textwidth}
\centering
\includegraphics[width=\textwidth]{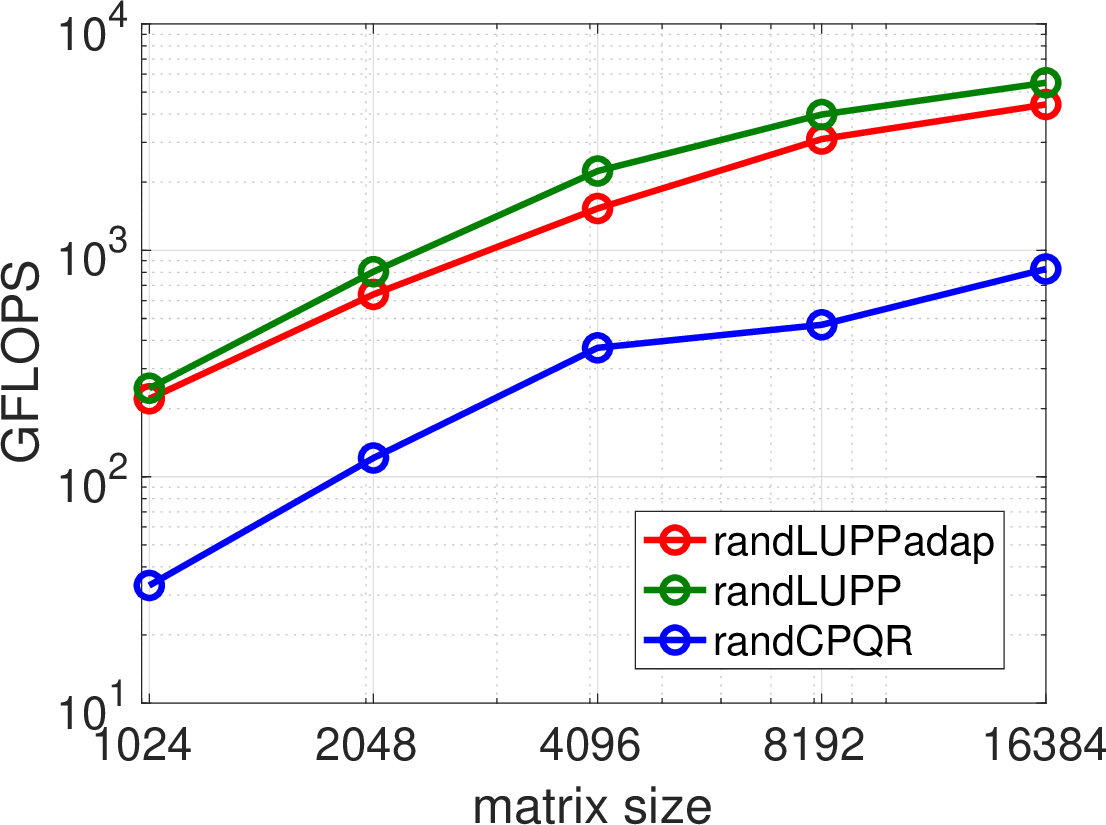}
\caption{Block size $b=128$}
\end{subfigure}
\begin{subfigure}{0.49\textwidth}
\centering
\includegraphics[width=\textwidth]{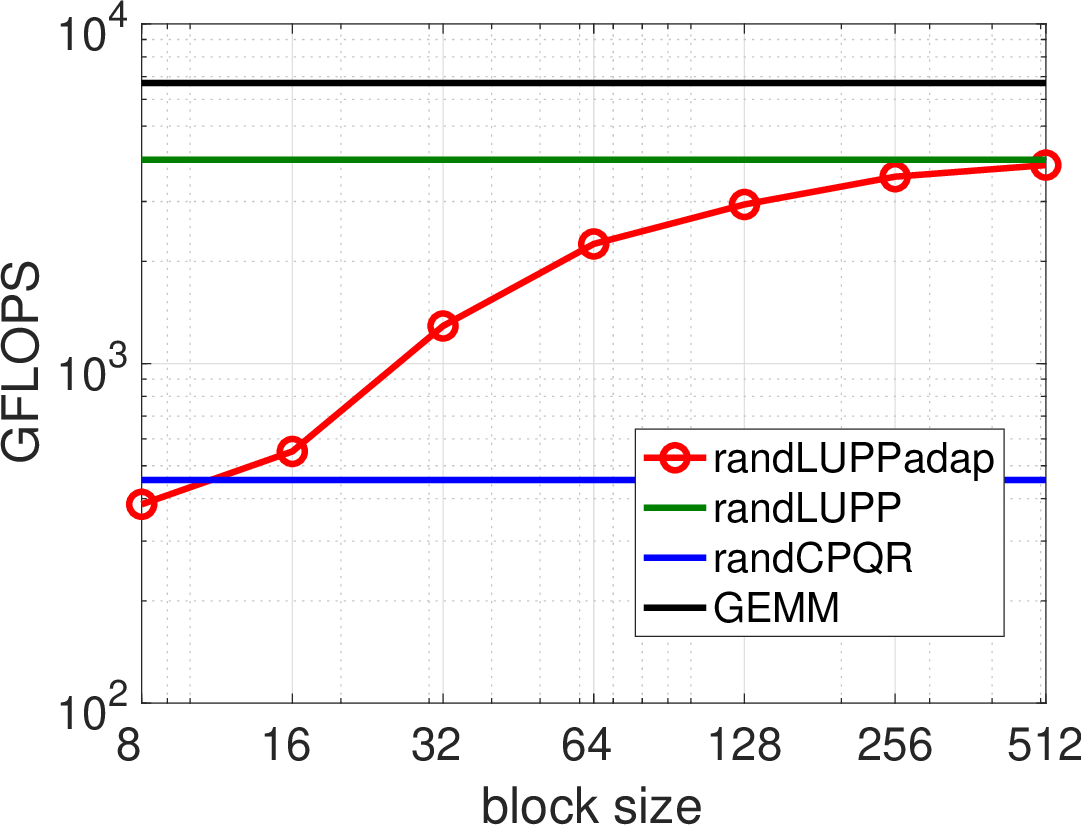}
\caption{Matrix size $n=8192$}
\end{subfigure}
\caption{GPU performance of \texttt{randLUPPadap} for increasing matrix sizes and for increasing block sizes, respectively. The test matrix is \textbf{Fast Decay} from Section~\ref{sec:acccomp}. The prescribed tolerance $\tau=10^{-8}$ corresponds to a numerical rank of approximately $n/2$, where $n$ is the matrix size. The computed numerical rank from \texttt{randLUPPadap} is used as an input for the two reference methods \texttt{randLUPP} and \texttt{randCPQR}. In (b), \texttt{GEMM} stands for multiplying two $n \times n$ matrices.
}
\label{f:nb}
\end{figure}

\begin{figure}
\centering
\begin{subfigure}{0.49\textwidth}
\centering
\includegraphics[width=\textwidth]{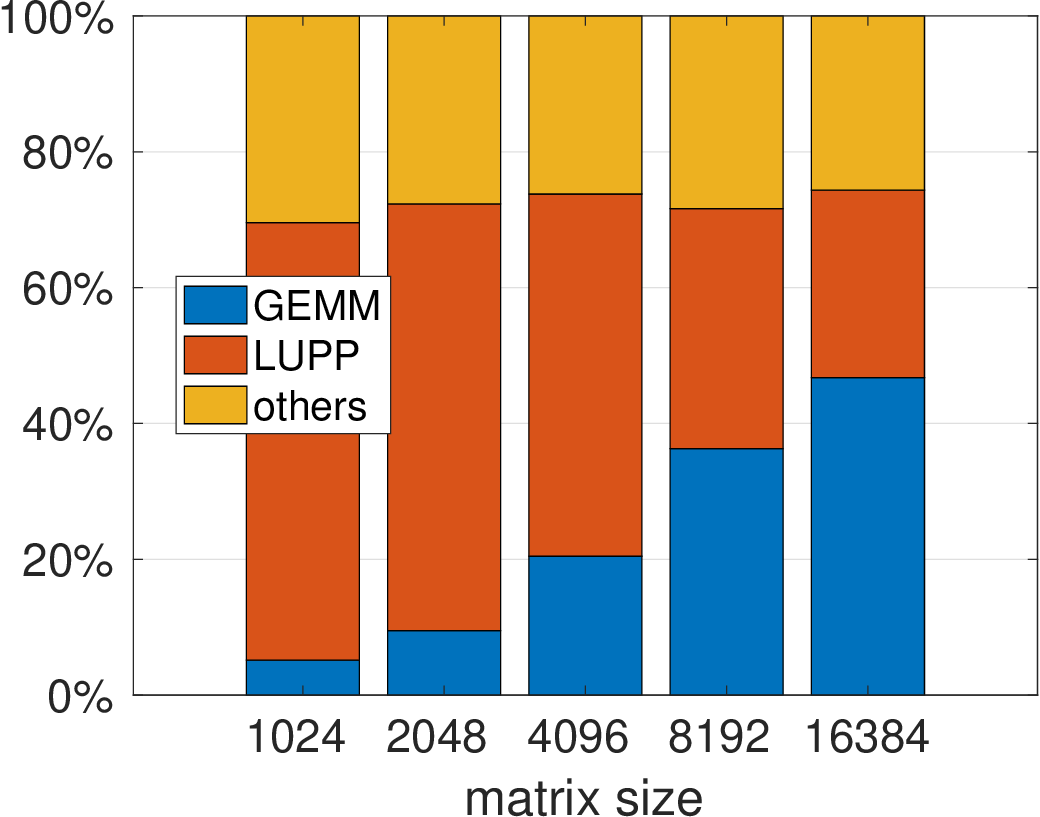}
\caption{Block size $b=128$}
\end{subfigure}
\begin{subfigure}{0.49\textwidth}
\centering
\includegraphics[width=\textwidth]{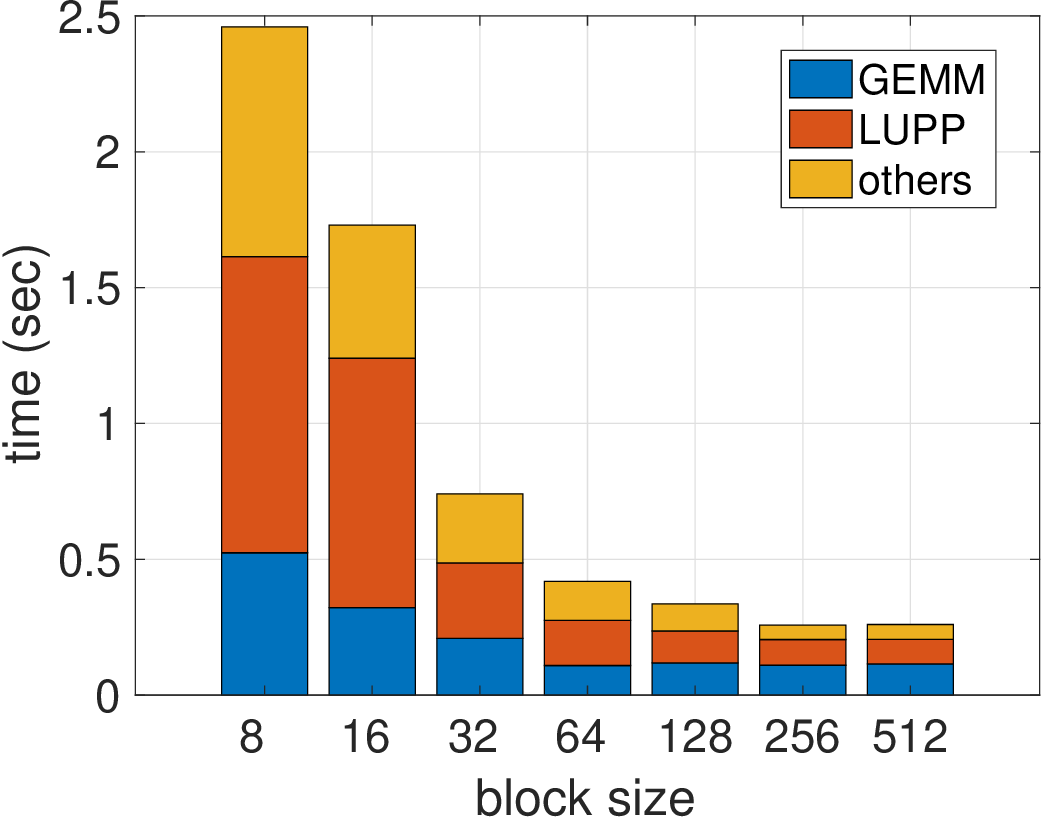}
\caption{Matrix size $n=8192$}
\end{subfigure}
\caption{
Profile of \texttt{randLUPPadap} for the GPU experiments in Figure~\ref{f:nb}. As shown in Algorithm~\ref{alg:IRMS}, matrix-matrix multiplication (GEMM) and LUPP are the most computationally intensive building blocks.
}
\label{f:prof}
\end{figure}

Figure~\ref{f:time} shows the running time of all three methods  on the GPU (solid lines) and on the CPU (dashed lines). On the GPU the new method \texttt{randLUPPadap} was almost as fast as \texttt{randLUPP}  and was approximately $6 \times$ faster than \texttt{randCPQR}  on average. 

The adaptive algorithm \texttt{randLUPPadap} is mathematically equivalent to \texttt{randLUPP} but works with small matrix blocks of size $b$. Therefore, the performance of \texttt{randLUPPadap} is upper bounded by that of  \texttt{randLUPP}, which is not far from that of matrix-matrix multiplication (GEMM).  The gap between \texttt{randLUPPadap} and \texttt{randLUPP}, however, can be small on the GPU  if the block size $b$ is sufficiently large, as shown in Figure~\ref{f:nb}.

Figure~\ref{f:nb} also shows that the new method  \texttt{randLUPPadap} can leverage the GPU better than  \texttt{randCPQR}  as long as the block size $b$ is not too small. The reason, as explained in Section~\ref{sec:intro} (see also Figure~\ref{f:luqr}), is that the computation of parallel {CPQR} is much less efficient than the computation of parallel {LUPP}.

Finally, Figure~\ref{f:prof} shows the runtime profile of  \texttt{randLUPPadap} for increasing input matrix sizes $n$ and increasing block sizes $b$, respectively. When the matrix size $n$ or the block size $b$ increases, all matrix computations in \texttt{randLUPPadap}  become more efficient (see Figure~\ref{f:prof} (b)), especially  the computation of parallel {LUPP}. This is why the runtime percentage of parallel {LUPP} decreases in Figure~\ref{f:prof} (a) when the matrix size $n$ increases.

\section{Conclusions and Future Work}
\label{sec:concl}

This paper introduces a randomized algorithm to construct an interpolative or CUR decomposition of an input matrix $\ma \in \F^{m \times n}$ adaptively using a prescribed error tolerance and sampling block size.
The algorithm leverages a cost-efficient method of error estimation using Schur complements with independent random samples.
Numerical experiments that investigate both accuracy and runtime performance demonstrate the algorithm's competitive degree of accuracy and improved efficiency over state-of-the-art methods.
Software is also provided for the adaptive parallelized algorithm for both GPU and CPU implementations.

Future avenues of exploration may include (1) investigating other cost-effective error estimation methods to compute the threshold metric; (2) applying the new method to construct the so-called rank-structured matrices~\cite{martinsson2019fast} arising from solvers for elliptic partial differential equations; and (3) extending the new method to compute the low-rank approximations of a batch of matrices, which frequently arises in computational and data science.

\section*{Acknowledgments}
KP was supported in part by the Peter O'Donnell Jr. Postdoctoral Fellowship at the Oden Institute.
CC was supported in part by startup funds of North Carolina State University.
YD was supported in part by the Office of Naval Research N00014-18-1-2354, NSF DMS-1952735, DOE ASCR DE-SC0022251, UT Austin Graduate School Summer Fellowship, and NYU Courant Instructorship.

\section{Appendix: Numerical Investigations of the LUPP Growth Factor}
\label{sec:appendix}

In this section, we record some of our findings regarding the asymptotic upper bound $\frac{k}{4 \log k}$ on the growth factor in LUPP for Haar-distributed matrices found in \cite{Higham21}.
We hypothesize that this upper bound holds for a larger class of random matrix distributions, supporting the recent conjecture in \cite{dong2023robust} and the analysis of sub-Gaussian random matrix distributions in \cite{saibaba2023randomized}.

As input matrices $\ma$, we used those in Section~\ref{sec:acccomp} (\textbf{Fast Decay}, \textbf{Kahan}, \textbf{large}, and \textbf{MNIST}), as well as the matrix presented by Chan in \cite{chan1987} exhibiting adversarial behavior with respect to the growth factor in LUPP:
\begin{align}
    \textbf{Chan}: \ma = \bmat{1 &  &  &   &  \\
    -1 & 1 &  &  &   \\
    -1 & -1 & 1 &  &    \\
    \vdots & \vdots &  & \ddots &    \\
    -1 & -1 & \cdots & -1  & 1}; \ \ m,n = 5000.
\end{align}

We used Gaussian, subsampled randomized trigonometric transforms (SRTT), and sparse sign (SS) random matrices to sketch each input matrix as in \ref{sec:acccomp}. 
We then computed $\mU_r$ as in Remark~\ref{sec:remark1}, using $p=10$ for oversampling.

\begin{figure}[h]
\centering
\begin{subfigure}{0.49\textwidth}
\centering
\includegraphics[width=\textwidth]{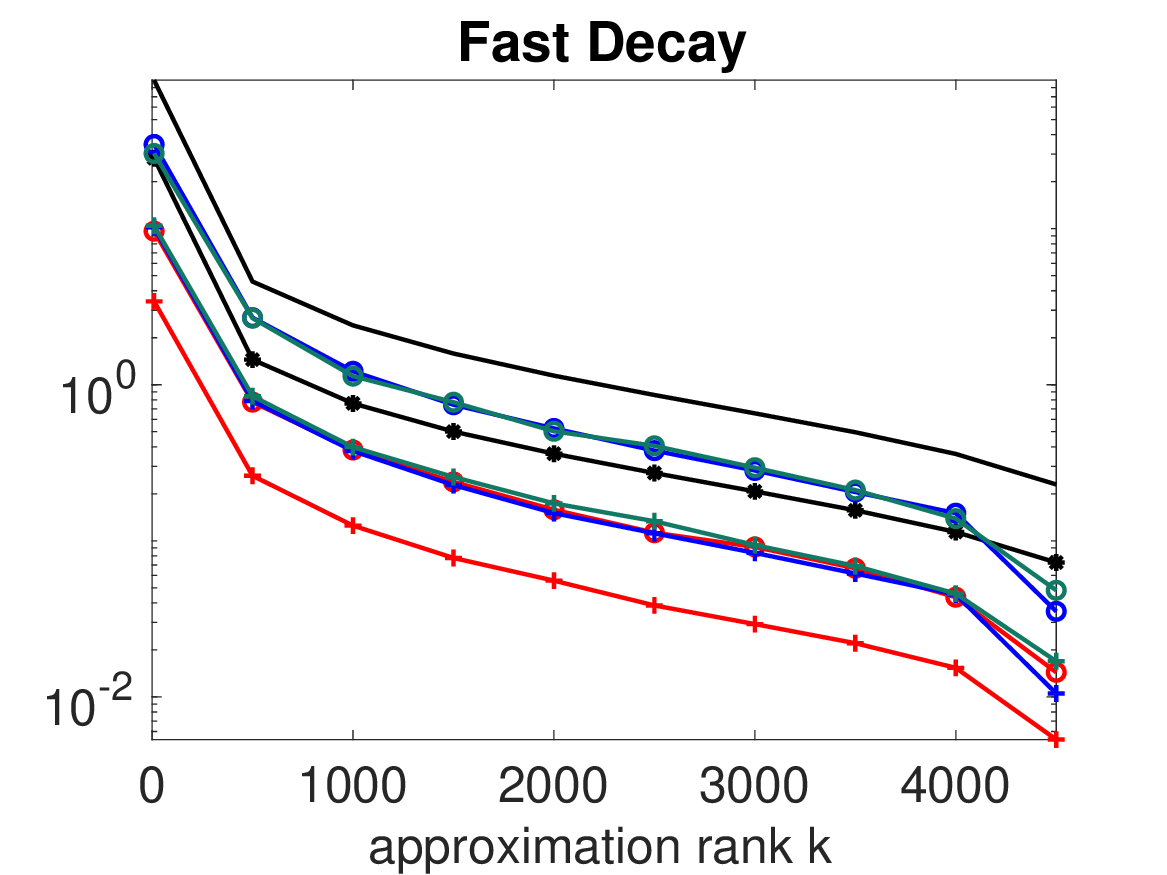}
\caption{}
\end{subfigure} 
\hfill
\begin{subfigure}{0.29\textwidth}
\includegraphics[width=\textwidth]{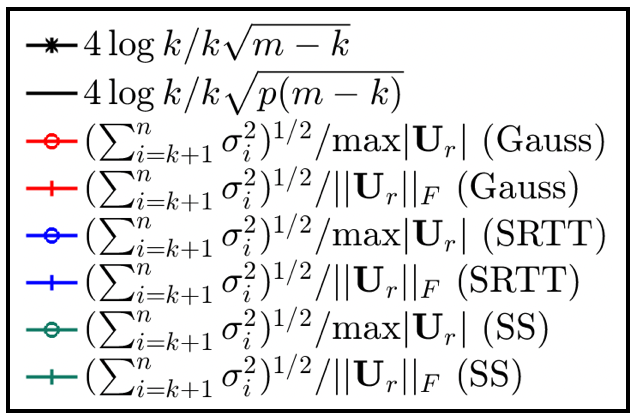}
\vspace{0mm}
\end{subfigure} 
\begin{minipage}{0.19\textwidth}
    \textcolor{white}{nothing}
\end{minipage}\\
\vspace{3mm}
\begin{subfigure}{0.49\textwidth}
\centering
\includegraphics[width=\textwidth]{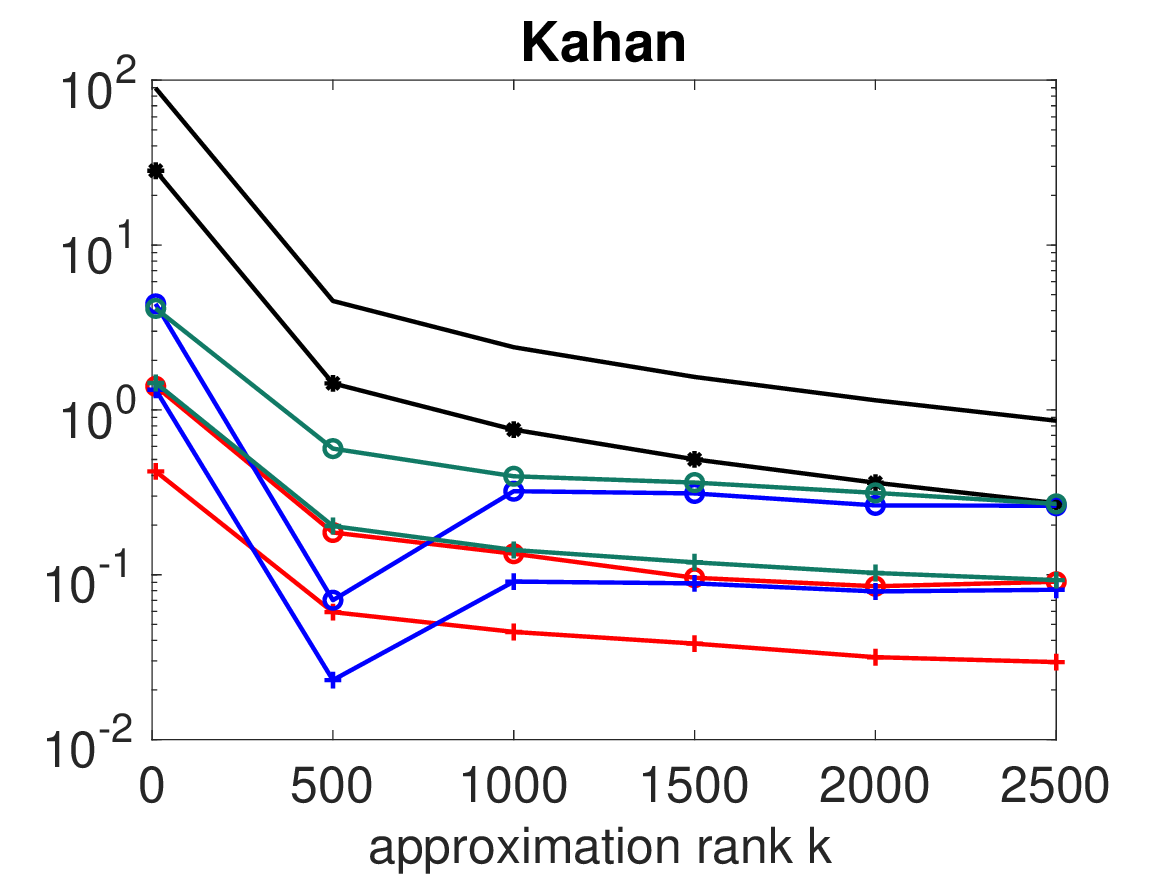}
\caption{}
\end{subfigure}
\begin{subfigure}{0.49\textwidth}
\centering
\includegraphics[width=\textwidth]{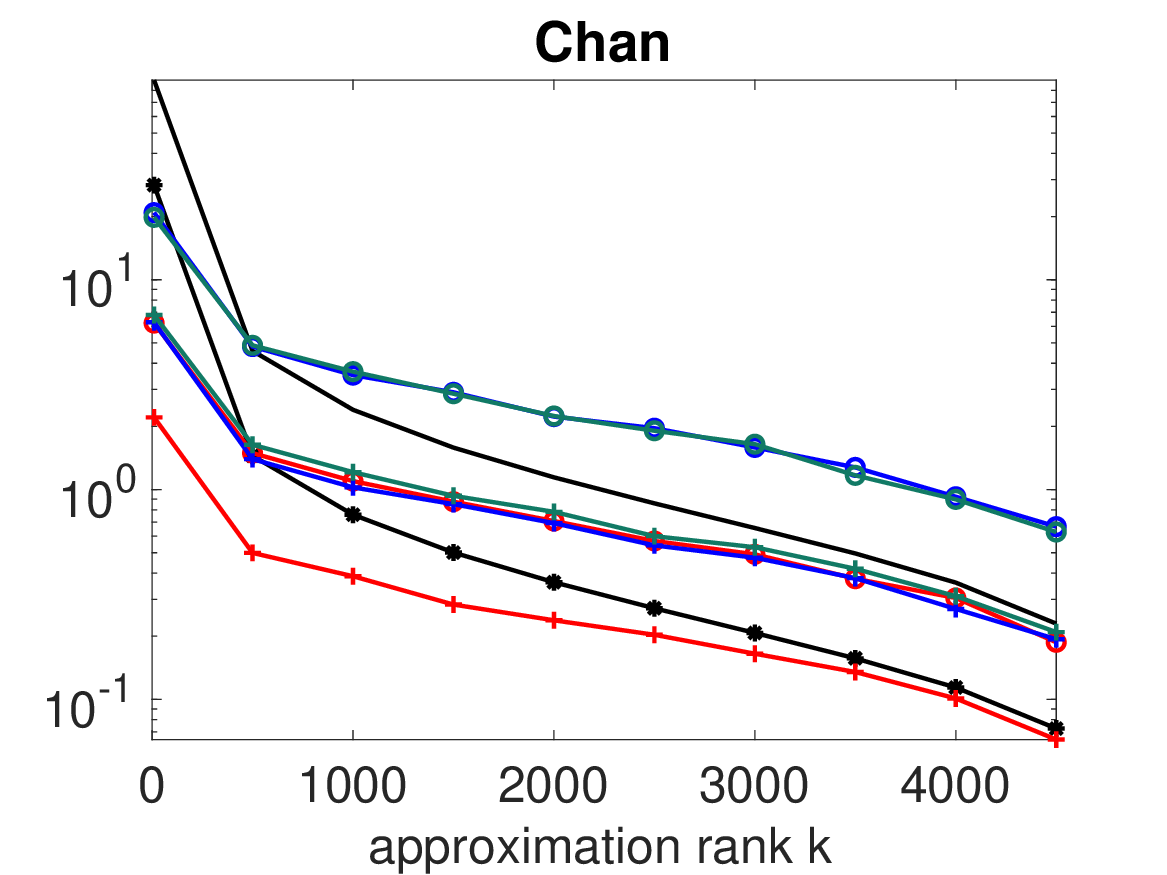}
\caption{}
\end{subfigure}
\begin{subfigure}{0.49\textwidth}
    \centering
    \includegraphics[width=\textwidth]{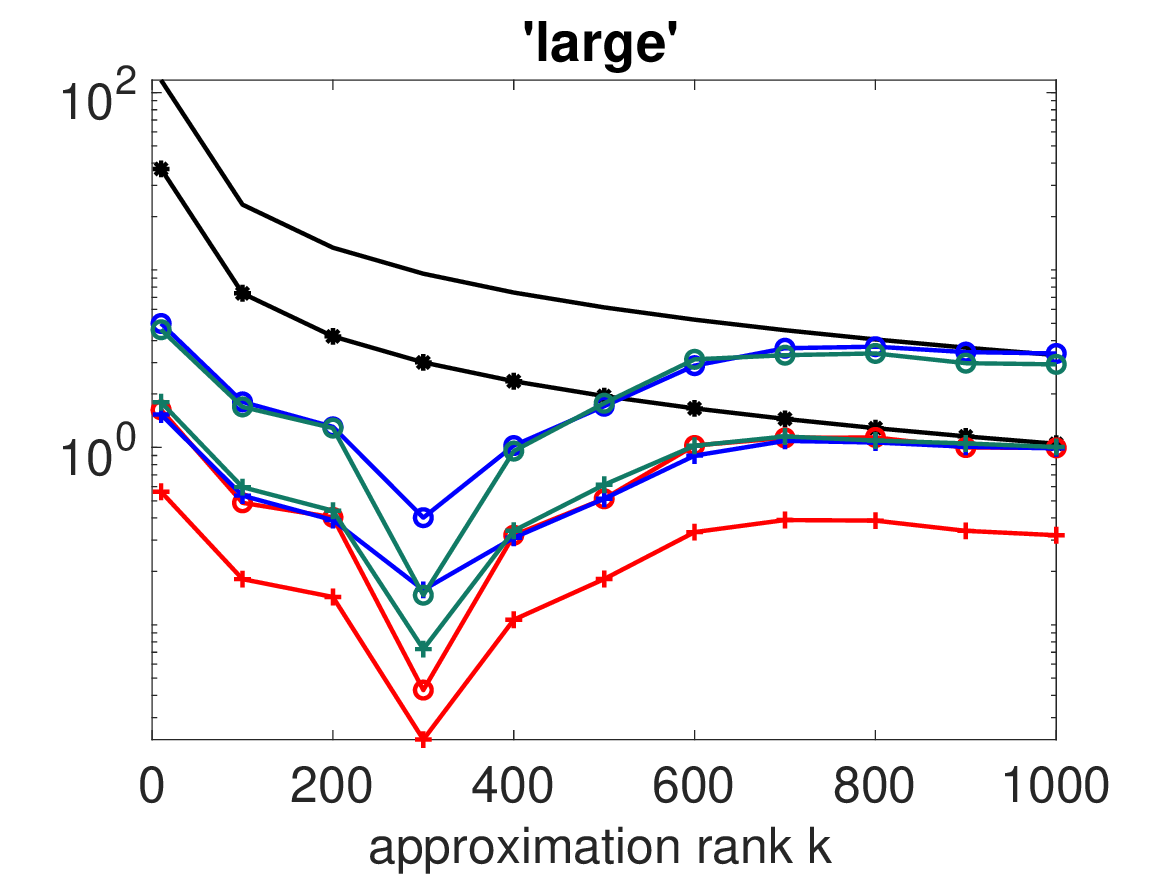}
\end{subfigure}
\begin{subfigure}{0.49\textwidth}
    \centering
    \includegraphics[width=\textwidth]{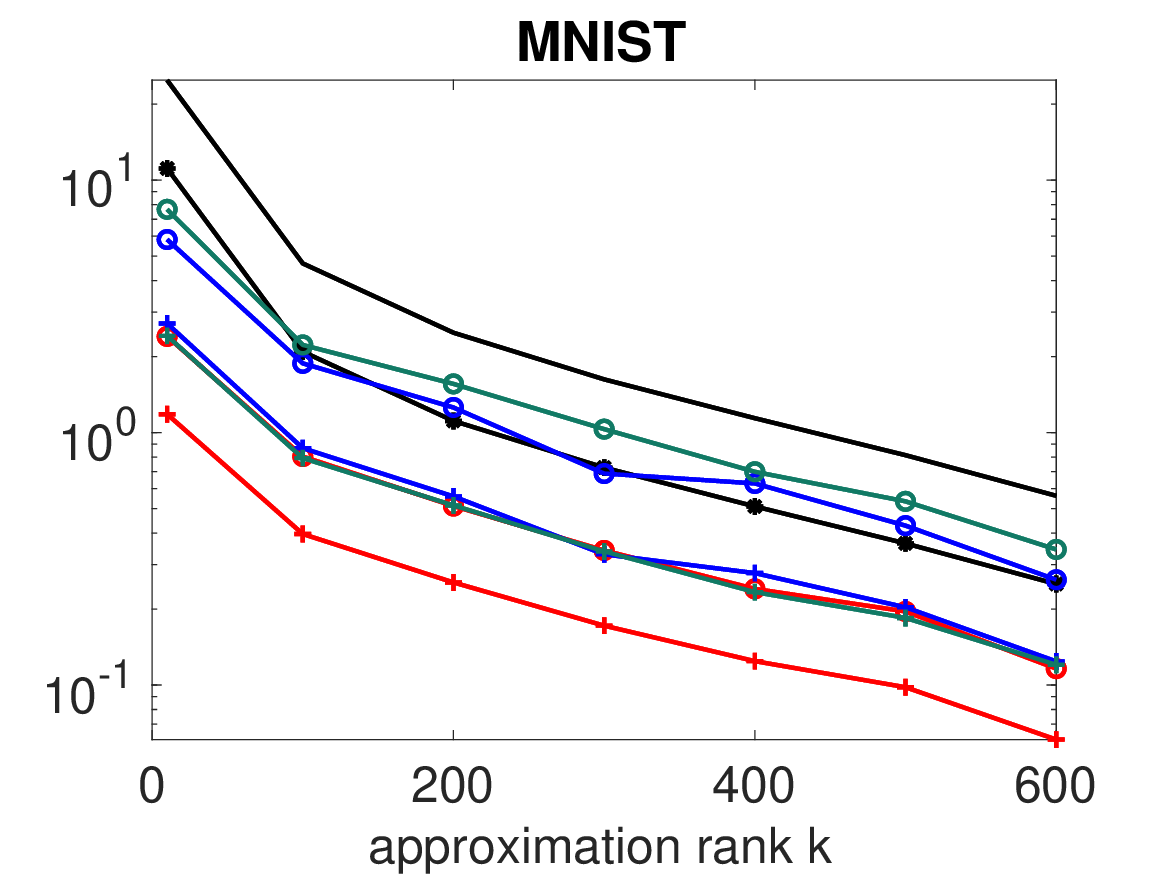}
\end{subfigure}
\caption{Ratios of the best rank-$k$ approximation errors to the computed means of $\max|\mU_r|$ and $\|\mU_r\|_F$. We compare these against the reciprocal of the asymptotic growth factor found in \cite{Higham21}, scaled by $\sqrt{p(m-k)}$ or $\sqrt{m-k}$.
}
\label{fig:ratios_synth}
\end{figure}

To study the asymptotic behavior, we computed the ratio of the best rank-$k$ approximation error to $\|\mU_r\|_F$ and $\max|\mU_r|$:
\begin{align}
    \frac{(\sum_{i=k+1}^n \sigma_i^2)^{1/2}}{\|\mU_r\|_F} \ \hspace{3mm} \mbox{and} \ \hspace{3mm} \frac{(\sum_{i=k+1}^n \sigma_i^2)^{1/2}}{\max|\mU_r|}.
\end{align}
In Figure~\ref{fig:ratios_synth}, we plot each of these values against $\frac{4 \log k}{k}$ scaled by $\sqrt{m-k}$ (for Gaussian) and $\sqrt{p(m-k)}$ (for SRTT, SS).
In all test cases, the reciprocal of the asymptotic upper bound on the growth factor yields an upper bound on the ratio of the SVD error to $\|\mU_r\|_F$, demonstrated by the plots with $+$ markers.

Moreover, with the exception of the adversarial Chan matrix (the likes of which are rarely encountered in practice \cite{chan1987}), the ratios of the SVD errors to the largest-magnitude elements of $\mU_r$ appear to converge to $\frac{4 \log k}{k}$. 
This behavior suggests that this asymptotic upper bound in combination with $\mU_r$ can be used to closely and reliably approximate the spectrum of the input matrix, for a ``rank-revealing'' LUPP (relying upon randomization).
More generally, it also suggests that this upper bound on the asymptotic growth factor holds more broadly than for Haar-distributed random matrices.

\bibliographystyle{abbrv}
\bibliography{ref}

\end{document}